\documentclass{amsart}
\usepackage{newlfont}
\allowdisplaybreaks

\let \SS=\S

\let\Min=\min
\let\Max=\max

        \newtheorem{theorem}{Theorem}
          \newtheorem{lemma}{Lemma}
      
     \newtheorem{definition}{Definition}
        \newtheorem{problem}{Problem}

\makeatletter%
\def\labelenumi{(\@roman\c@enumi)}%
\def\theenumi{(\@roman\c@enumi)}%
\def\labelenumii{(\@alph\c@enumii)}%
\def\theenumii{\@alph\c@enumii}%
\def\labelenumiii{(\@arabic\c@enumiii)}%
\def\theenumiii{(\@arabic\c@enumiii)}%
\def\p@enumiii{\theenumi\theenumii}%
\def\wideA#1{\@mathmeasure\z@\textstyle{#1}%
\ifdim\wd\z@>\tw@ em\mathaccent "055C{#1}%
\else\mathaccent "0364{#1}\fi}%
\def\wideB#1{\@mathmeasure\z@\textstyle{#1}%
\ifdim\wd\z@>\tw@ em\mathaccent "055E{#1}%
\else\mathaccent "0366{#1}\fi}%
\makeatother%
\def\makecs#1#2{\makecsX {#1}#2,.}%
\def\makecsX#1#2#3.{\onecs{#1}{#2}%
  \ifx#3,\let\next\eatit\else\let\next\makecsX\fi\next{#1}#3.}%
\def\onecs#1#2{\expandafter\gdef\csname #2\endcsname%
  {{\csname #1\endcsname {#2}}}}%
\def\eatit#1#2.{\relax}%
\makecs{}{abcdefghijklmnopqrstuvwxyzABCDEFGHIJKLMNOPQRSTUVWXYZ}%
\makecs{mathbb}{Q}%
\def\Zn{\mathbb{Z}}%
\def\blank{\phantom{o}}%
\def\con#1{\setbox13\hbox{$#1$}\ifdim\wd13<1em\breve{#1}\else{\left(#1\right)}%
  \breve{\ }\fi}%
\def\rp{{\hts;\hts}}%

\def\id{{1\kern-.08em\raise1.3ex\hbox{\normalshape\rm,}\kern.08em}}%
\def\di{{0\kern-.04em\raise1.3ex\hbox{\normalshape\rm,}\kern.04em}}%
\def\halfthinspace{\relax\ifmmode\mskip.5\thinmuskip\relax%
  \else\kern.8888em\fi}%
\let\hts=\halfthinspace%
\def\Id{\mathsf{Id}\hts}%
\def\Di{\mathsf{Di}\hts} \def\ifff{\ \ \text{iff}\ \ }%
\def\conv#1{\setbox13\hbox{$#1$}\ifdim\wd13<1.1em{#1}^{-1}%
  \else{\left(#1\right)}^{-1}\fi}%
\def\ll#1{L^I_{#1}}%
\def\cl#1{R^I_{#1}}%
\def\ll#1{L_{#1}}%
\def\cl#1{R_{#1}}%
\def\rmin{{\sim}}%

\def\RM{\text{\normalshape\sf RM}}%
\def\BM{\text{\normalshape\sf BM}}%
\def\CL{\text{\normalshape\sf CL}}%
\def\KR{\text{\normalshape\sf KR}}%
\def\RR{\text{\normalshape\sf R}}%
\def\EE{\text{\normalshape\sf E}}%
\def\MM{\text{\normalshape\sf M}_0}%

\def\<{\langle}%
\def\>{\rangle}%
\def\ie{{\it i.e.}\/}%
\def\min#1{\overline{#1}}%
\def\Su{\mathbf{S}}%
\def\sra{\mathfrak\S}%
\def\At{\mathcal{A}t}%
\def\Ch{\mathcal{C}}%
\def\KRM{K_\RM}%
\def\Ka{\mathfrak{K}}%
\def\alg#1#2{$#1_{#2}$}%
\def\SI#1{\S^\I_{#1}}%
\def\TI#1{{\T^\I_{#1}}}%
\def\TH#1{{\hat\T^\I_{#1}}}%
\def\Rm{\mathcal{R}}%
\def\3{\emptyset}%
\def\Re{\mathfrak{Re}} \def\Sb{\mathfrak{Sb}}%
\def\gc#1{\mathfrak{#1}}%
\def\at{{At}\hts}%
\def\Cm#1{\setbox13\hbox{$#1$}{\gc{Cm}({#1})}}%
\begin{document}
\title[{\sf Sugihara relation algebras}] {Relation algebras of
  Sugihara, Belnap, Meyer, and Church} \author{R.\ L.\ Kramer, R.\ D.\
  Maddux}
\address{Department of Mathematics\\
  396 Carver Hall\\ Iowa State University\\ Ames, IA 50011\\ U.S.A.}
\email{maddux@iastate.edu} \subjclass{Primary: 03G15, 03B47}
\keywords{relevance logic, Sugihara chains, Belnap's $M_0$, Meyer's
  RM84, Church's diamond, crystal lattice, relation algebra, proper
  relation algebra, representable relation algebra, dense relation,
  transitive relation, $R$-mingle, definitional reduct} \date{\today}
\begin{abstract}
  Algebras introduced by, or attributed to, Sugihara, Belnap, Meyer,
  and Church are representable as algebras of binary relations with
  set-theoretically defined operations.  They are definitional reducts
  or subreducts of proper relation algebras.  The representability of
  Sugihara matrices yields sound and complete set-theoretical
  semantics for $\RR$-mingle.
\end{abstract}
\maketitle
%\tableofcontents
%\listoffigures
%\listoftables
\section{Introduction}\label{sect1}
Sugihara's matrix, described by A.\ Anderson and N.\ Belnap \cite[pp.\
335--6]{MR0406756}, was introduced by T.\ Sugihara in 1955
\cite{Sugihara1955}. A smaller one, obtained by using only one element
per integer instead of two, is taken by Anderson and Belnap as ``the
Sugihara matrix''.  R.\ K.\ Meyer introduced finite Sugihara matrices
for his proof that they are complete for the Dunn-McCall logic
$\RR$-mingle, or $\RM$ \cite[\SS29.3.2]{MR0406756}.  Various algebras,
including all Sugihara matrices and perhaps others, are representable
as algebras of binary relations.  Their operations are defined
set-theoretically, and need not be specified by tables.  Since their
operations are definable in the similarity type of relation algebras,
they are definitional reducts or subreducts of proper relation
algebras.  This was proved already in \cite{MR2641636} for finite
Sugihara matrices of even cardinality. In this paper we extend this
result to all finite Sugihara matrices plus Sugihara's original
infinite matrix and two others described by Anderson and Belnap. We
also show it for Belnap's $\MM$ and for matrices of Meyer and Church.
These algebras may be represented by a list of relations on a set,
together with some operations on relations selected from Table
\ref{table1}.  We start with Belnap.
\begin{table}
  \begin{itemize}
  \item[$\bullet$] identity relation on $\U$,
    $\Id=\{\<\x,\x\>:\x\in\U\}$,
  \item[$\bullet$] diversity relation on $\U$, 
    $\Di=\{\<\x,\y\>:\x,\y\in\U,\,\x\neq\y\}$,
  \item[$\bullet$] universal relation on $\U$,
    $\U^2=\{\<\x,\y\>:\x,\y\in\U\}$,
  \item[$\bullet$] union, $\A\cup\B=\{\<\x,\y\>:\<\x,\y\>\in\A
    \text{ or } \<\x,\y\>\in\B\}$,
  \item[$\bullet$] intersection, $\A\cap\B=\{\<\x,\y\>:
    \<\x,\y\>\in\A \text{ and } \<\x,\y\>\in\B\}$,
  \item[$\bullet$] converse,
    $\conv\A=\{\<\x,\y\>:\<\y,\x\>\in\A\}$,
  \item[$\bullet$] complement,
    $\min\A=\{\<\x,\y\>:\x,\y\in\U,\,\<\x,\y\>\notin\A\}$,
  \item[$\bullet$] converse-complement,
    $\rmin\A=\{\<\x,\y\>:\x,\y\in\U,\,\<\y,\x\>\notin\A\}$,
  \item[$\bullet$] relative product, \newline
    $\A|\B=\{\<\x,\y\>:$ for some $\z\in\U$, $\<\x,\z\>\in\A$ and
    $\<\z,\y\>\in\B\}$,
  \item[$\bullet$] residual, \newline $\A\to\B=\{\<\x,\y\>:$ for
    all $\z\in\U$, if $\<\z,\x\>\in\A$ then $\<\z,\y\>\in\B\}$,
  \item[$\bullet$] relativized converse-complement, \newline
    $\rmin'\A=\{\<\x,\y\>:\<\y,\x\>\in\Di$ and
    $\<\y,\x\>\notin\A\}$,
  \item[$\bullet$] relativized relative product,
    $\A|'\B=\{\<\x,\y\>: \<\x,\y\>\in\Di$, and \newline for some
    $\z\in\U$, $\<\x,\z\>\in\A\cap\Di$ and $\<\z,\y\>\in\B\cap\Di\}$,
  \item[$\bullet$] relativized residual,
    $\A\to'\B=\{\<\x,\y\>:\big(\<\x,\y\>\in\Di$, and \newline for
    all $\z\in\U$, if $\<\z,\x\>\in\A\cap\Di$ and $\<\z,\y\>\in\Di$
    then $\<\z,\y\>\in\B\big)\}$.
  \end{itemize}
  \caption{Some relations on a set $U$ and some operations on relations.}
\label{table1}
\end{table}
\section{Belnap}\label{sect2}
\begin{figure}
  \begin{picture}(40,40)(0,0)%
    \put(20,30){\circle*{1}{$\Q^2$}}%
    \put(16.5,30){$+3$}%
    \put(40,20){\circle*{1}}%
    \put(41,19.5){$\geq$}%
    \put(36,19.5){$+2$}%
    \put(20,20){\circle*{1}{$\neq$}}%
    \put(16.5,20){$-0$}%
    \put(40,10){\circle*{1}}%
    \put(41,9.5){{$>$}}%
    \put(36,9.5){{$-2$}}%
    \put( 0,10){\circle*{1}}%
    \put(-2.5,9.5){$<$}%
    \put(1.5,9.5){$-1$}%
    \put(20,10){\circle*{1}}%
    \put(21,9){=}%
    \put(16.5,9){+0}%
    \put( 0,20){\circle*{1}}%
    \put(-2.5,19.5){$\leq$}%
    \put(1.5,19.5){$+1$}%
    \put(20, 0){\circle*{1}}%
    \put(21,-1){$\emptyset$}%
    \put(16.5,-1){$-3$}%
    \put(20, 0){\line(2,1){20}}%
    \put( 0,10){\line(2,1){20}}%
    \put(20,10){\line(2,1){20}}%
    \put( 0,20){\line(2,1){20}}%
    \put( 0,10){\line(2,-1){20}}%
    \put(20,30){\line(2,-1){20}}%
    \put(20,20){\line(2,-1){20}}%
    \put( 0,20){\line(2,-1){20}}%
    \put( 0,10){\line(0, 1){10}}%
    \put(20,20){\line(0, 1){10}}%
    \put(40,10){\line(0, 1){10}}%
    \put( 20,0){\line(0, 1){10}}%
\end{picture}
\caption{Lattice of relations in $M_0$.}
\label{fig1}
\end{figure}
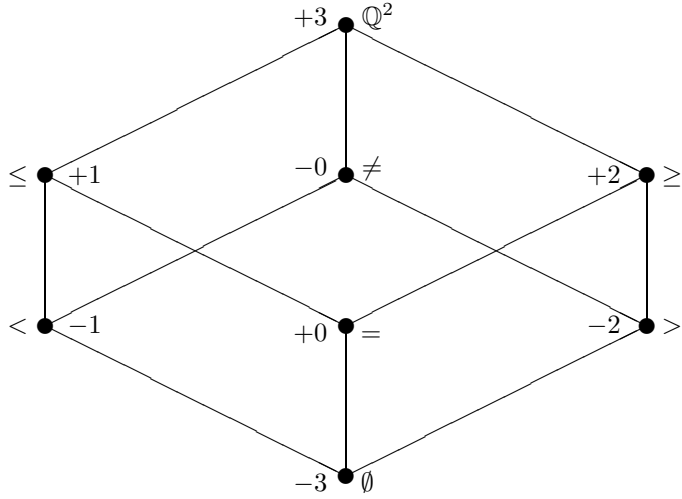

Belnap's $\MM$ was first introduced in 1960 \cite{MR0141590} by
matrices for binary operations $\lor$, $\land$, $\to$, $\rmin$, and
unary operations $N$ and $M$, on an eight-element set.  From the
matrices for $\land$ and $\lor$ it is apparent that the eight values
appearing in them, namely $-3$, $-2$, $-1$, $-0$, $+0$, $+1$, $+2$,
and $+3$ (the last four are the designated values), form a lattice
isomorphic to the lattice of subsets of the 3-element set
$\{-1,+0,-2\}$, with $+3$ at the top and $-3$ at the bottom, where
$\land$ and $\lor$ are interpreted as intersection and union.  This
observation does not occur in~\cite{MR0141590}, but in subsequent
literature $\MM$ is usually portrayed this way, by a Hasse diagram
along with tables for $\to$ and $\sim$. See, for example,
\cite[\SS18.4, \SS22.1.3]{MR0406756}, \cite[\SS\SS34.1--2]{MR1223997},
\cite[p.\,178]{MR728950}, or \cite[pp.\ 101--2]{MR3728341}.  It is
described in \cite[p.\ 117]{MR2536403}, \cite{Maddux2007}, and
\cite[Theorem~4.1]{MR2641636} as an algebra
\begin{align}\label{belnap}
  \MM&=\<\M_0,\cup,\,\cap,\,\to,\,\rmin\>,&
  \M_0&=\{\emptyset,\,<,\,>,\,=,\,\neq,\,\leq,\,\geq,\,\Q^2\},
\end{align}
whose universe $\M_0$ consists of eight binary relations on the
rational numbers $\Q$: the empty relation $\emptyset$, the less-than
relation $<$, the greater-than relation $>$, the identity relation
$=$, the diversity relation $\neq$, less-than-or-equal $\leq$,
greater-than-or-equal $\geq$, and the universal relation $\Q^2$. These
eight relations are the unions of subsets of $\{<,\,>,\,=\}$.  The
Hasse diagram for $\M_0$ is shown in Figure \ref{fig1}.
The operations of $\MM$ are union $\cup$, intersection $\cap$,
residuation $\to$, and converse-complementation $\rmin$, defined in
Table \ref{table1} with $\U=\Q$.  The logic called $\BM$
\cite[p.\,128]{MR3728341} is defined by an explicit finite
axiomatization. By \cite[Theorem 9.8.6]{MR3728341} and its corollary,
$\MM$ is characteristic for the logic $\BM$.  Because it has a single
finite characteristic structure, $\BM$ is a complete decidable logic.
The universe $\M_0$ of $\MM$ is also closed under complementation
$\min\blank$, conversion $\conv{}$, and relative multiplication $|$,
and contains the empty relation $\emptyset$, universal relation
$\Q^2$, and identity relation $\Id$ on $\Q$.  Therefore $\M_0$ is the
universe of an algebra
\begin{align*}
  \mathfrak\M_0&=\<\M_0,\cup,\,\cap,\,\min\blank,\,\emptyset,\,\Q^2,
  \,|,\,\conv{},\,\Id\>.
\end{align*}
We refer to $\mathfrak\M_0$ as {\bf Belnap's relation algebra}.  It is
a proper relation algebra on the set of rational numbers.  Proper
relation algebras were first defined in \cite{JonssonTarski1948},
\cite[Definition 4.23]{MR45086}, and \cite[\SS2]{MR0037278}.
\begin{definition}\label{proper}
  For any equivalence relation $\E$, let
  \begin{align*}
    \Sb(\E)&=\<\wp(\E),\,\cup,\,\cap,\,\min\blank\,,\,
    \emptyset,\,\E,|,\conv{},\,\Id\>,
  \end{align*}
  where $\wp(\E)$ is the set of all subsets of $\E$, $\cup$ is union,
  $\cap$ is intersection, $\min\blank$ is complementation with respect
  to $\E$, $|$ is relative multiplication, $\conv{}$ is conversion,
  $\emptyset$ is the empty relation, and $\Id = \{ \<\u,\u\> :
  \<\u,\u\> \in \E\}$ is the identity relation on the field of $\E$.
  $\Sb(\E)$ is the {\bf algebra of subrelations of $\E$}.  A {\bf
    proper relation algebra} is any subalgebra of the algebra of
  subrelations of an equivalence relation.  An algebra is {\bf
    representable} if it is isomorphic to a proper relation algebra.
  For any set $\U$, let
  \begin{align*}
    \Re(\U)&=\Sb(\U^2).
  \end{align*}
  $\Re(\U)$ is the {\bf algebra of relations on $\U$}.  A {\bf proper
    relation algebra on $\U$} is any subalgebra of $\Re(\U)$.
\end{definition}
$\Re(\U)$ is the prototypical example of a relation algebra. Tarski's
original axioms \cite{MR5280} were chosen because they are true in
$\Re(\U)$. Tarski's axiom XII implies simplicity, that is, any algebra
satisfying XII has no non-trivial homomorphic images. This axiom was
later dropped in \cite{MR0043763, JonssonTarski1948} so that all
proper relation algebras would satisfy the axioms for relation
algebras.  It was noticed very early that Belnap's $\MM$ is a lattice
with additional operations. What required nearly half a century after
its introduction in 1960 was the realization, first mentioned in 2007
\cite{Maddux2007}, that $\MM$ is a definitional reduct of
$\mathfrak\M_0$ (see Definition \ref{reducts}).  The proper relation
algebra $\mathfrak\M_0$ was already known to Lyndon in 1956
\cite{MR79570}. In footnote 13, p.\ 307, Lyndon says,
\begin{quote}
  ``Every relation algebra without zero divisors that is of order not
  exceeding 8 (there are 13 such) is commutative and isomorphic to a
  complex algebra of either the additive rationals or a cyclic group
  of order not exceeding 13.''
\end{quote}
In the numbering system of \cite{MR2269199}, the 13 relation algebras
without zero divisors of order 8 or less are algebra \alg11 of order
2, algebras \alg12 and \alg22 of order 4, and the ten algebras \alg13,
\alg23, \alg33, \alg17, \alg27, \alg37, \alg47, \alg57, \alg67, and
\alg77 of order 8.  Algebras \alg11, \alg12, \alg23, \alg17, and
\alg57 are isomorphic to proper relation algebras on sets of size 1,
2, 3, 4, and 5, respectively, but on no larger or smaller sets.
Algebras \alg22, \alg27, \alg47, \alg67, and \alg77 are isomorphic to
proper relation algebras on sets of size at least 3, 6, 9, 8, and 9,
respectively, and also on sets of all larger sizes.  Algebra \alg47
shows up in \SS\ref{sect8}, where it is called $\mathfrak{Ch}$. It has
Church's diamond as a definitional reduct.  Algebra \alg37 is
isomorphic to proper relation algebras on all sets of even cardinality
6 or larger.  Algebra \alg33 also shows up in \SS\ref{sect8}, where it
is called $\mathfrak{Rm}$.  It has Meyer's algebra RM84 as a
definitional reduct. Algebra \alg33 is isomorphic to proper relation
algebras on sets of cardinality 7, and 9 or more, but not on sets of
size 8 \cite[Theorem 4.2]{MR1334290}.  The representation on 7
elements appears in \SS\ref{sect8}.  Finally, algebra \alg13 is
isomorphic to Belnap's relation algebra $\mathfrak\M_0$. It is
isomorphic to proper relation algebras on sets of every infinite
cardinality. All the representations of these algebras on the smallest
possible sets are unique \cite{MR1334290}.  The representation on the
rationals $\Q$ mentioned by Lyndon is unique because of Cantor's
theorem on the categoricity of dense linear orderings without
endpoints on countable sets. It is the one relation algebra mentioned
by Lyndon that requires the ``additive rationals''.

Starting with \cite{Allen1981, Allen1984, Allen1983}, an extensive
literature developed in the 1980s in which $\mathfrak\M_0$ is known as
the {\bf Point Algebra}, because among its eight relations are the
three ways that two points on the rational number line can be related
to each other: either they are equal ($=$), or the first point is to
the left of the second point ($<$), or to the right ($>$). The Point
Algebra and similar algebras based on the relationships that hold
between various combinations of points and regions are widely used in
computer science for spatial and temporal reasoning, and for
constraint satisfaction problems. Consult \cite{MR1369206}, where
references to some of the early work can be found.  More recent papers
that explicitly mention the Point Algebra include \cite{MR2235660,
  MR2225706, MR3873385, MR2565927, MR2722825, MR1909067, MR2015214,
  MR3466223, MR1658912, MR2158576, MR2790871, MR1872076, MR3944722,
  MR3089976, MR1927631}.
\section{Sugihara}\label{sect3}
A {\bf lattice} is an algebra $\<\S,\lor,\land\>$ with binary
operations $\lor$ and $\land$ that are associative, commutative, and
idempotent, such that the absorption laws $\A\land(\A\lor\B) = \A =
\A\lor(\A\land\B)$ hold. A lattice is a {\bf chain} if $\A\land\B$ is
always either $\A$ or $\B$, \ie, the ordering $\leq$ is linear, where
$\A\leq\B\ifff\A\land\B=\A$.
\begin{definition}\label{def1}
  $\Su=\<\S,\lor,\land,\to,\rmin\>$ is a {\bf Sugihara chain} if
  $\<\S,\lor,\land\>$ is a chain, $\rmin$ is an involution that
  reverses the ordering, \ie,
  \begin{equation*}
    \rmin\rmin\A=\A,\qquad\qquad\A\leq\B\ifff\rmin\B\leq\rmin\A,
  \end{equation*}
  and $\A\to\B=\rmin\A\lor\B$ if $\A\leq\B$, otherwise
  $\A\to\B=\rmin\A\land\B$. An element $\A\in\S$ is said to be {\bf
    designated} if $\rmin\A\leq\A$.
\end{definition}
Ten examples of Sugihara chains residing in Belnap's $\M_0$ are shown
in Table \ref{ten}, specified by their relations and operations.
\begin{table}
  \begin{align*}
    &\text{Universe}&&\text{Operations}\\
    &\{<,\leq\}, \quad \{>,\geq\}, \quad \{\emptyset,\Q^2\},
    &&\cup,\,\cap,\,\to,\,\rmin,\\
    &\{\emptyset,<,\leq,\Q^2\}, \quad \{\emptyset,>,\geq,\Q^2\},%
    &&\cup,\,\cap,\,\to,\,\rmin,\\
    &\{<\}, \quad \{>\}, \quad \{\emptyset,\neq\},
    &&\cup,\,\cap,\,\to',\,\rmin',\\
    &\{\emptyset,<,\neq\}, \quad \{\emptyset,>,\neq\},
    &&\cup,\,\cap,\,\to',\,\rmin'.
  \end{align*}
  \caption{Ten examples of Sugihara chains in $M_0$}\label{ten}
\end{table}
The chains $\{<,\leq\}$ and $\{>,\geq\}$ appear in Belnap's original
proof \cite[p.\,145]{MR0141590} of the variable-sharing property for
the logic $\EE$ of Anderson-Belnap \cite[\SS21.1]{MR0406756}, which
says that if $\A\to\B$ is a theorem of $\EE$ then $\A$ and $\B$ share
at least one propositional variable. The same proof applies to the
logic $\RR$.  Axioms (R1)--(R13) for $\RR$ are shown in Table
\ref{table2}; see \cite[\SS27.1.1]{MR0406756} or \cite[pp.\
xxiii--xxvi]{MR1223997}.  The axioms of $\RR$ are valid in $\MM$, and
the rules of deduction preserve validity, so $\A\to\B$ is not a
theorem of $\RR$ whenever $\A$ and $\B$ share no variable.

For every finite cardinality there is exactly one Sugihara chain
having that cardinality.  Because $\rmin$ is order-reversing, finite
Sugihara chains of odd cardinality must have an element that is a
fixed point for $\rmin$, the one in the middle. Such an element would
be assigned as a truth value to a formula that is equivalent to its
own negation.  Sugihara chains without fixed points under negation are
called ``normal'' by Meyer~\cite[p.\ 400]{MR0406756}, so odd Sugihara
chains are not normal.  Sugihara chains with even cardinality have no
elements fixed by $\rmin$, which interchanges the top and bottom
halves while reversing their order. Sugihara chains with even
cardinality were used by Meyer \cite[p.\ 413, Corollary
3.1]{MR0406756} to prove that the theorems of $\RM$ are exactly those
formulas that are valid in all finite Sugihara chains.  His result was
used to prove \cite[Theorem 6.2(iii)]{MR2641636}, which says that a
formula is a theorem of $\RM$ if and only if it is valid in every
finite algebra in $\KRM$ (see Definition \ref{KRM}).

Infinite Sugihara chains are not determined by cardinality alone; see
\cite{MR3526497}. In what Anderson and Belnap ``have accordingly come
to think of \dots\ as \emph{the} Sugihara matrix'' \cite[p.\
337]{MR0406756}, the universe is the set
$\Zn^*=\{\n\colon0\neq\n\in\Zn\}$ of non-zero integers and
$\rmin(i)=-i$ for every non-zero $i\in\Zn^*$.  This Sugihara chain was
named $\Su_{\Zn^*}$ by Meyer~\cite[p.\ 414]{MR0406756}, who proved
that the theorems of $\RM$ are exactly those formulas valid in
$\Su_{\Zn^*}$ \cite[p.\ 414, Corollary 3.5]{MR0406756}.  Having
described $\Su_{\Zn^*}$, Anderson and Belnap suggest, ``Or one might
insert 0 between $-1$ and $+1$, counting it designated''
\cite[p.\,337]{MR0406756}. The resulting chain was called $\Su_\Zn$ by
Meyer \cite[p.\ 414]{MR0406756}. It has a fixed point for $\rmin$,
namely $0=\rmin0$.  No such fixed point occurs in the original chain
of Sugihara \cite{Sugihara1955}.  In this chain, the ordering is
isomorphic to two copies of the integers, one after the other, so we
call it $\Su_{\Zn+\Zn}$ (with $+$ denoting ordinal addition).  In more
detail, the elements are $s_i$ and $t_j$ for integers
$i,j\in\mathbb\Z$, and the ordering is defined by $s_i<s_j$ and
$t_i<t_j$ whenever $i<j$, and $s_i<t_j$ for any $i$ and $j$.  The
Sugihara chains $\Su_{\Zn^*}$, $\Su_\Zn$, and $\Su_{\Zn+\Zn}$ are
countably infinite but not isomorphic. We turn now to the construction
of proper relation algebras that have these and all finite Sugihara
matrices as definitional reducts.
\section{Definition of $\mathcal{S}_I$}\label{sect4}
For an arbitrary index set $\I\subseteq\Zn$ of integers,
$\mathcal\S_\I$ is a set of relations on ${}^\Zn\Q$, where ${}^\Zn\Q$
be the set of functions $\q:\Zn\to\Q$ that map the integers $\Zn$ to
the rationals $\Q$.  By Theorem \ref{lemma1} in the next section,
$\mathcal\S_\I$ is the universe of a proper relation algebra called
$\sra_\I$.  When $\I=\{0\}$, Belnap's $\MM$ is a definitional reduct
of $\sra_{\{0\}}$ (which is isomorphic to Belnap's relation algebra),
and, when $\I=\Zn$, Sugihara's original matrix $\Su_{\Zn+\Zn}$ is a
definitional subreduct (but not a definitional reduct) of $\sra_\Zn$.
\begin{definition}\label{defU} Let $\I\subseteq\Zn$.
  \begin{enumerate}
  \item If $\q\in{}^\Zn\Q$, we say that $\q$ is {\bf eventually zero}
    if there exists some integer $\n\in\Zn$ such that $\q_\i=0$ for
    every integer $\i>\n$.
  \item Let $\U_\I$ be the set of functions, called {\bf sequences},
    that map $\Zn$ to $\Q$, are eventually zero, and are non-zero only
    on $\I$:
    \begin{align*}
      \U_\I=\{\q\colon\q\in{}^\Zn\Q,\,
      (\exists_{\n\in\Zn})(\forall_{\i>\n})(\q_\i=0),\,
      (\forall_{\i\in\Zn})(\i\notin\I\Rightarrow\q_\i=0)\}.
    \end{align*}
  \item Define the identity and diversity relations
    \begin{align*}
      \Id_\I&=\{\<\q,\q\>\colon\q\in\U_\I\},\qquad
      \Di_\I=\{\<\q,\r\>\colon\q,\r\in\U_\I,\,\q\neq\r\},\\
      \intertext{and, for every $\n\in\Zn$,}
      \ll\n&=\{\<\q,\r\>\colon\text{$\q,\r\in\U_\I$, $\q_\n<\r_\n$,
        and $\q_\i=\r_\i$ whenever $\n<\i$}\},\\
      \cl\n&=\{\<\q,\r\>\colon\text{$\q,\r\in\U_\I$, $\q_\n>\r_\n$,
        and $\q_\i=\r_\i$ whenever $\n<\i$}\}.
    \end{align*}
  \item Define a set of relations and its set of unions
    \begin{align*}
      \At_\I&=\{\Id_\I\}\cup\bigcup_{\n\in\I}\{\ll\n,\,\cl\n\},&
      \mathcal\S_\I&=\left\{\bigcup\mathcal\X\colon
        \mathcal\X\subseteq\At_\I\right\}.
    \end{align*}
\end{enumerate}
\end{definition}
By Theorem \ref{lemma1} below, $\At_\I$ is the set of atoms of a
complete atomic proper relation algebra called $\sra_\I$.  Let
$\I=\emptyset$. Then $\U_\emptyset$ is the set consisting of just the
one sequence $\q=\<\cdots,0,0,0,\cdots\>$ that is always zero.  Notice
that for every $\n\in\Zn$, $\ll\n=\emptyset$ iff $\cl\n=\emptyset$ iff
$\n\notin\I$.  Hence $\At_\emptyset=\{\Id_\emptyset\}=\{\<\q,\q\>\}$
and $\mathcal\S_\emptyset=\{\emptyset,\{\q\}\}$, so $\sra_\emptyset$
is isomorphic to the proper relation algebra $\Re(\{\q\})$.

Suppose $\I=\{0\}$. In this case $\U_{\{0\}}$ is the set of
$\Zn$-indexed sequences of rational numbers having $0$ everywhere
except possibly at index $0$.  There is a bijection between
$\U_{\{0\}}$ and $\Q$ that maps $\q\in\U_{\{0\}}$ to $\q_0\in\Q$.
Setting $\I=\{0\}$ in Definition \ref{defU} gives
\begin{align*}
  \At_{\{0\}}&=\{\Id_{\{0\}},\ll0,\,\cl0\},\\
  \mathcal\S_{\{0\}} &=\left\{\emptyset,\,\Id_{\{0\}},\,\ll0,\,\cl0,\,
    \Id_{\{0\}}\cup\ll0,\,\Id_{\{0\}}\cup\cl0,\,
    \ll0\cup\cl0,\,(\U_{\{0\}})^2\right\}.
\end{align*}
For every relation $\X\subseteq(\U_{\{0\}})^2$, let
$\f(\X)=\{\<\q_0,\r_0\>\colon\<\q,\r\>\in\X\}$. Applying $\f$ to the
relations in $\At_{\{0\}}$ produces the relations in $\M_0$:
\begin{align*}
  \f(\emptyset)&=\emptyset,\\
  \f(\Id_{\{0\}})&=\{\<\x,\y\>\colon\x,\y\in\Q,\,\x=\y\},\\
  \f(\ll0)&=\{\<\x,\y\>\colon\x,\y\in\Q,\,\x<\y\},\\
  \f(\cl0)&=\{\<\x,\y\>\colon\x,\y\in\Q,\,\x>\y\},\\
  \f(\Id_{\{0\}}\cup\ll0)&=\{\<\x,\y\>\colon\x,\y\in\Q,\,\x\leq\y\},\\
  \f(\Id_{\{0\}}\cup\cl0)&=\{\<\x,\y\>\colon\x,\y\in\Q,\,\x\geq\y\},\\
  \f(\ll0\cup\cl0)&=\{\<\x,\y\>\colon\x,\y\in\Q,\,\x\neq\y\},\\
  \f(\U_{\{0\}})&=\Q^2.
\end{align*}
In fact, $\f$ is an isomorphism from $\sra_{\{0\}}$ to Belnap's
relation algebra, so $\sra_{\{0\}}$ is also called ``Belnap's relation
algebra''. It contains copies of the ten Sugihara chains in Table
\ref{ten}.

When $\I=\Zn$ the set $\mathcal\S_\Zn$ contains far more than is
needed for the original Sugihara chain. $\At_\Zn$ has countably many
relations, so the cardinality of $\mathcal\S_\Zn$ is same as that of
the real numbers. By Theorem \ref{lemma1} in the next section,
$\mathcal\S_\Zn$ is the universe of a relation algebra $\sra_\Zn$,
called {\bf Sugihara's relation algebra}.  By Theorem \ref{iso} in the
section after that, $\sra_\Zn$ contains countable chains isomorphic to
$\Su_{\Zn+\Zn}$.
\section{Structure of $\mathfrak{S}_I$}\label{sect5}
In this section we show $\mathcal\S_\I$ is the universe of the
complete atomic proper relation algebra $\sra_\I$.  First we review
Definition \ref{defU} for easy reference.  For every set of integers
$\I\subseteq\Zn$, $\U_\I$ is the set of functions from $\Zn$ to $\Q$
that are eventually zero and non-zero only on $\I$, $\At_\I =
\{\Id_\I\} \cup \bigcup_{\n\in\I} \{\ll\n,\,\cl\n\}$, and
$\mathcal\S_\I = \{\bigcup \mathcal\X\colon \mathcal\X \subseteq
\At_\I\}$, where $\Id_\I = \{ \<\q,\q\> \colon \q\in\U_\I \}$, and for
every $\n\in\I$,
\begin{align*}
  &\ll\n=\{\<\q,\r\>\colon\text{$\q,\r\in\U_\I$, $\q_\n<\r_\n$, and
    $\q_\i=\r_\i$ whenever $\n<\i$}\},\\
  &\cl\n=\{\<\q,\r\>\colon\text{$\q,\r\in\U_\I$, $\q_\n>\r_\n$, and
    $\q_\i=\r_\i$ whenever $\n<\i$}\}.
\end{align*} 
\begin{theorem}\label{lemma1}
  $\At_\I$ is a partition of $(\U_\I)^2$. $\At_\I$ is the set of atoms
  of the complete atomic proper relation algebra
  \begin{align*}
    \sra_\I&=\<\mathcal\S_\I,
    \cup,\cap,\min\blank,\emptyset,(\U_\I)^2,\,|,\conv{},\Id_\I\>.
  \end{align*}
\end{theorem}
\proof To see that the relations in $\At_\I$ are pairwise disjoint and
their union is $(\U_\I)^2$, note that any two sequences
$\q,\r\in\U_\I$ are either equal everywhere (are in the identity
relation $\Id_\I$), or differ somewhere, in which case there is a
\emph{largest} integer $\n$ where they differ, since they are both
eventually zero.  The pair $\<\q,\r\>$ cannot be in $\ll\m$ or $\cl\m$
if $\n<\m$ since $\q$ and $\r$ agree at every such $\m$ by the choice
of $\n$, and $\<\q,\r\>$ cannot be in $\ll\m$ or $\cl\m$ if $\n>\m$
since $\q$ and $\r$ differ somewhere larger than $\m$, namely at $\n$.
Since $\q$ and $\r$ differ at $\n$, one of them is not zero at $\n$,
so $\n\in\I$.  Since the ordering of the rationals is linear, either
$\q_\n<\r_\n$ or $\q_\n>\r_\n$ but not both, so the pair $\<\q,\r\>$
must be in the relation $\ll\n$ or $\cl\n$ but not both.  Therefore we
have a disjoint union:
\begin{align*}
  (\U_\I)^2&=\bigcup_{\X\in\At_\I}\X.
\end{align*}
Since the relations in $\At_\I$ form a partition of $(\U_\I)^2$, the
unions of arbitrary subsets of $\At_\I$ form a complete atomic Boolean
algebra whose set of atoms is $\At_\I$.  Thus $\mathcal\S_\I$ is
closed under union, intersection, complementation, and contains
$\emptyset$, $(\U_\I)^2$, and $\Id_\I$.  It remains to be verified
that $\mathcal\S_\I$ is closed under conversion and relative
multiplication. From their definitions it follows that $\cl\n$ and
$\ll\n$ are converses of each other. The converse of the identity
relation is itself. Conversion distributes over arbitrary unions of
relations. We therefore have the following rules.  For all $\n\in\I$
and $\mathcal\X\subseteq\At_\I$,
\begin{align}\label{rule-3}
  \conv{\ll\n}&=\cl\n,& \conv\Id_\I&=\Id_\I,&
  \conv{\bigcup\mathcal\X}&=\bigcup_{\X\in\mathcal\X}\conv\X.
\end{align}
From \eqref{rule-3} it follows that
$\{\bigcup\mathcal\X\colon\mathcal\X\subseteq\At_I\}$ is closed under
conversion.  For closure under relative multiplication, we reason as
follows.  The relative product of two unions of sets of atoms is, by
distributivity, the union of the relative products of the atoms in the
two sets. More exactly, if $\mathcal\X,\mathcal\Y\subseteq\At_\I$ then
\begin{align*}
  \bigcup\mathcal\X|\bigcup\mathcal\Y&=
  \bigcup\{\X|\Y\colon\X\in\mathcal\X,\,\Y\in\mathcal\Y\}.
\end{align*}
The relative product is again a union of atoms if the relative product
of any two atoms is a union of atoms. As we will see, the relative
product of any two atoms is an atom in every case except the relative
product of a diversity atom and its converse, in which case the
relative product is the union of the identity relation and all the
diversity atoms with smaller index; see \eqref{rule0}.

Assume $\q$, $\r$, and $\s$ are distinct sequences in $\U_\I$.  Each
sequence must differ from the other two, so the cardinality of
$\{\q_\n,\r_\n,\s_\n\}$ cannot be $1$ for every $\n\in\Zn$. On the
other hand, since $\q$, $\r$, and $\s$ are all eventually zero, the
number of elements of $\{\q_\n,\r_\n,\s_\n\}$ will eventually be
constantly $1$, since $\{\q_\n,\r_\n,\s_\n\}=\{0\}$ whenever $\n$ is
large enough.  Hence there is an integer $\n$ at which
$\{\q_\n,\r_\n,\s_\n\}$ contains either exactly three or exactly two
elements (hence $\n\in\I$ because they can't all be zero) and
$|\{\q_\i,\r_\i,\s_\i\}|=1$ for all $\i>\n$ (hence $\q$, $\r$, and
$\s$ all agree beyond $\n$).  Although any two of $\q$, $\r$, and $\s$
are equal beyond $\n$, any pair of them could also agree beyond an
integer strictly smaller than $\n$.  In the first case, when
$\{\q_\n,\r_\n,\s_\n\}$ has exactly three elements, those elements
must form a chain under the dense linear ordering $<$ on the
rationals, and, since $\q$, $\r$, and $\s$ all agree beyond $\n$, we
may choose $\x,\y,\z$ so that $\{\x,\y,\z\}=\{\q,\r,\s\}$ and
$\x\mathrel{\ll\n}\y\mathrel{\ll\n}\z$ (and $\x\mathrel{\ll\n}\z$).
This is listed as case \eqref{case1} below.  If there are exactly two
elements in $\{\q_\n,\r_\n,\s_\n\}$, then one of them differs from the
other two, and the other two coincide.  Therefore, for some $\x,\y,\z$
such that $\{\x,\y,\z\}=\{\q,\r,\s\}$, we have $\x_\n\neq\y_\n=\z_\n$
and $1=|\{\x_\i,\y_\i,\z_\i\}|$ for every $\i>\n$. If
$\x_\n<\y_\n=\z_\n$ then $\x\mathrel{\ll\n}\y$ and
$\x\mathrel{\ll\n}\z$, while if $\x_\n>\y_\n=\z_\n$ then
$\x\mathrel{\cl\n}\y$ and $\x\mathrel{\cl\n}\z$.  Now $\y$ and $\z$
are distinct, but they agree beyond $\n$ and also agree at $\n$.
Hence they disagree at some $\j<\n$, and agree beyond $\j$, in which
case $\y\mathrel{\ll\j}\z$ or $\y\mathrel{\cl\j}\z$.  We may assume
$\x,\y,\z$ were chosen so that $\y\mathrel{\ll\j}\z$.  This yields the
remaining two cases \eqref{case2} and \eqref{case3}.  Thus, given any
three distinct $\q,\r,\s\in\U_\I$, there are $\x,\y,\z\in\U_\I$ and
$\n\in\I$ such that $\{\x,\y,\z\}=\{\q,\r,\s\}$,
$|\{\q_\n,\r_\n,\s_\n\}|>1$, $|\{\q_\i,\r_\i,\s_\i\}|=1$ for all
$\i>\n$, and one of these three cases holds:
\begin{align}
  \label{case1}
  &\x\mathrel{\ll\n}\y\mathrel{\ll\n}\z\text{ and }
  \x\mathrel{\ll\n}\z,\\
  \label{case2}
  &\x\mathrel{\ll\n}\y\mathrel{\ll\j}\z\text{ and }
  \x\mathrel{\ll\n}\z\text{ for all }\j<\n,\\
  \label{case3}
  &\x\mathrel{\cl\n}\y\mathrel{\ll\j}\z\text{ and }
  \x\mathrel{\cl\n}\z\text{ for all }\j<\n.
\end{align}
From the fact that these are the only possible cases, we will be able
to deduce the rules for computing relative products of pairs of
relations in $\At_\I$.  First we consider the relative products with
the identity relation.
\begin{align}\label{rule-2}
  \Id_\I|\Id_\I&=\Id_\I,&\ll\n|\Id_\I=\Id_\I|\ll\n&=\ll\n,
  &\cl\n|\Id_\I=\Id_\I|\cl\n&=\cl\n.
\end{align}
We will only prove $\Id_\I|\ll\n=\ll\n$.  The other equations have
similar proofs. Assume $\<\q,\r\>\in\Id_\I|\ll\n$. Then there is some
$\s$ such that $\<\q,\s\>\in\Id_\I$ and $\<\s,\r\>\in\ll\n$. The
latter two statements tell us that $\q=\s$ and $\s\mathrel{\ll\n}\r$,
from which we conclude $\q\mathrel{\ll\n}\r$ by the fact that equal
objects have the same properties, hence $\<\q,\r\>\in\ll\n$, showing
that $\Id_\I|\ll\n\subseteq\ll\n$.  For the opposite inclusion, we
assume $\<\q,\r\>\in\ll\n$ and note that by choosing $\s=\q$ we get
$\<\q,\s\>\in\Id_\I$ and $\<\s,\r\>\in\ll\n$, hence $\<\q,\r\> \in
\Id_\I|\ll\n$.  Thus $\Id_\I\subseteq\Id_\I|\ll\n$.  Combining this
with $\Id_\I|\ll\n\subseteq\Id_\I$, we obtain the desired equality.

Next we introduce notation for special relations in $\mathcal\S_\I$
that arise from relative products of diversity atoms.  For any
$\n,\m\in\I$ let
\begin{align}
  \label{notation1}
  \ll{[\n,\m]}&=\bigcup\{\ll\k\colon{\n\leq\k\leq\m},\,\k\in\I\},&
  \ll{(-\infty,\n]}&=\bigcup\{\ll\k\colon\n\geq\k\in\I\},\\
  \label{notation2}
  \ll{[\n,\infty)}&=\bigcup\{\ll\k\colon\n\leq\k\in\I\},&
  \ll{(-\infty,\infty)}&=\bigcup\{\ll\k\colon\k\in\I\}.
\end{align}
Note that $\ll{[\n,\m]}=\emptyset$ if $\n>\m$, and $\ll{[\n,\n]} =
\ll\n$.  The same notation is used with converses (change $\L$ to $\R$
in the equations above). The rules \eqref{rule-3} imply
\begin{align*}
  \conv{\ll{[\n,\m]}}&=\cl{[\n,\m]},&
  \conv{\ll{(\infty,\m]}}&=\cl{(\infty,\m]},\\
  \conv{\ll{[\n,\infty)}}&=\cl{[\n,\infty)},&
  \conv{\ll{(-\infty,\infty)}}&=\cl{(-\infty,\infty)}.
\end{align*}
The relative product of diversity atoms $\ll\m,\cl\n\in\At_\I$ can be
computed according to four basic rules.  Rule \eqref{rule-1} says that
the relative product of a diversity atom with itself is itself.  Rule
\eqref{rule0} says that the relative product of a diversity atom with
its converse is the union of the identity relation and all the
diversity atoms having equal or smaller index. Rules \eqref{rule1} and
\eqref{rule2} say that the relative product of two diversity atoms
with distinct indices $\n$ and $\m$ is the one with the larger index.
\begin{align}\label{rule-1}
  \ll\n|\ll\n&=\ll\n,\qquad\cl\n|\cl\n=\cl\n,\\
  \label{rule0}
  \cl\n|\ll\n&=\ll\n|\cl\n=\Id_\I
  \cup\ll{(-\infty,\n]}\cup\cl{(-\infty,\n]},\\
  \label{rule1}
  \ll\m|\ll\n&=\ll\n|\ll\m=\cl\m|\ll\n=
  \ll\n|\cl\m=\ll\n
  \quad\text{if $\m<\n$},\\
  \label{rule2}
  \cl\m|\cl\n&=\cl\n|\cl\m=\cl\n|\ll\m= \ll\m|\cl\n=\cl\n %
  \quad\text{if $\m<\n$.}
\end{align}
To prove \eqref{rule-1}, assume $\<\q,\r\>\in\ll\n|\ll\n$, so there is
some $\s\in\U_\I$ such that $\<\q,\s\>\in\ll\n$ and
$\<\s,\r\>\in\ll\n$.  It follows that $\q_n<\s_\n$, $\s_\n<\r_\n$, and
$\q$, $\r$, and $\s$ all agree beyond $\n$.  We also have $\q_n<\r_\n$
by the transitivity of the ordering $<$ on $\Q$, so
$\<\q,\r\>\in\ll\n$. This shows $\ll\n|\ll\n\subseteq\ll\n$.  For the
opposite inclusion, assume $\<\q,\r\>\in\ll\n$. Then $\q_\n<\r_\n$,
and $\q$ and $\r$ agree beyond $\n$.  By the density of $<$, we may
choose $\s\in\U_\I$ so that $\s$ agrees with $\q$ and $\r$ beyond $\n$
and has some value $\s_\n$ between $\q_\n$ and $\r_\n$ (such as the
average of $\q_\n$ and $\r_\n$) so that $\q_\n<\s_\n<\r_\n$. The
values of $\s$ on arguments in $\I$ and smaller than $\n$ are
arbitrary.  This completes the proof of the first equation in
\eqref{rule-1}. The second equation has a similar proof.

To prove \eqref{rule0}, assume $\<\q,\r\>\in\cl\n|\ll\n$. If $q=\r$
then $\<\q,\r\>$ is in $\Id_\I$, one of the relations in the union on
the right side of \eqref{rule0}, as desired, so assume $\q\neq\r$.  By
the definition of $|$ there is some $\s\in\U_\I$ such that
$\<\q,\s\>\in\cl\n$, $\<\s,\r\>\in\ll\n$, and $\q$, $\r$, and $\s$
agree beyond $\n$.  From $\<\s,\q\>\in\ll\n$, $\<\s,\r\>\in\ll\n$, and
$\q\neq\r$ we conclude that we are in case \eqref{case1} or
\eqref{case2} with $\s=\x$ and $\{\q,\r\}=\{\y,\z\}$.  Since
$\q\neq\r$ and $\{\q,\r\}=\{\y,\z\}$, $\<\q,\r\>$ is in some diversity
atom whose index must be either $\n$, as in case \eqref{case1}, or
some smaller integer $\j<\n$, which occurs in case \eqref{case2}. In
either case, depending on how $\q$ and $\r$ match up with $\y$ and
$\z$, we have $\<\q,\r\>\in\ll{(-\infty,\n]}\cup\cl{(-\infty,\n]}$.
The pair $\<\q,r\>$ thus belongs to one of the relations on the right,
proving
\begin{align*}
  \cl\n|\ll\n\subseteq\Id_\I\cup\ll{(-\infty,\n]}\cup\cl{(-\infty,\n]}.
\end{align*}
For the converse, assume $\<\q,\r\>\in\Id_\I\cup\ll\m\cup\cl\m$ and
$\m\leq\n$.  We will find $\s\in\U_\I$ such that $\<\q,\s\>\in\cl\n$
and $\<\s,\r\>\in\ll\n$.  Since $\q$ and $\r$ agree beyond $\m$, they
agree beyond $\n$ as well.  Choose values for $\s\in\U_\I$ so that
$\s$ agrees with $\q$ and $\r$ beyond $\n$. The values of $\s$ at
arguments that are in $\I$ and smaller than $\n$ may be anything.  At
$\n$, choose a rational $\s_\n$ that is strictly smaller than both
$\q_\n$ and $\r_n$, such as $\s_\n=\Min(\q_\n,\r_n)-1$.  Here we are
using the fact that the ordering of the rationals does not have any
endpoints.  From $\s_\n<\q_\n$, $\s_\n<\r_n$, and the agreement of
$\q$, $\r$, and $\s$ beyond $\n$ we get $\<\s,\q\>\in\ll\n$ and
$\<\s,\r\>\in\ll\n$, hence $\<\q,\s\>\in\cl\n$, so $\<\q,\r\> \in
\cl\n|\ll\n$. This shows $\cl\n|\ll\n \supseteq \Id_\I \cup
\ll{(-\infty,\n]} \cup \cl{(-\infty,\n]}$, completing the proof of one
of the two equations in \eqref{rule0}. The other equation may be
proved similarly.

The proofs for \eqref{rule1} and \eqref{rule2} are somewhat simpler.
We first show that if $\m<\n$ then
$(\ll\m\cup\cl\m)|\ll\n\subseteq\ll\n$.  Assume
$\<\q,\r\>\in(\ll\m\cup\cl\m)|\ll\n$.  Then there must exist some
$\s\in\U_\I$ such that $\<\q,\s\>\in\ll\m\cup\cl\m$ and
$\<\s,\r\>\in\ll\n$.  It follows from $\<\q,\s\>\in\ll\m\cup\cl\m$
that $\q_\m\neq\s_\m$ and $\q$ and $\s$ agree beyond $\m$. Since
$\m<\n$, this tells us that $\q_\n=\s_\n$ and $\q$ and $\s$ agree
beyond $\n$.  From $\<\s,\r\>\in\ll\n$ we know $\s_\n<\r_\n$ and $\s$
and $\r$ agree beyond $\n$.  We conclude that $\q_\n=\s_\n<\r_\n$ and
$\q$, $\r$, and $\s$ agree beyond $\n$, hence $\<\q,\r\>\in\ll\n$.
Assume $\<\q,\r\>\in\ll\n$ and $\m<\n$.  Then $\q_\n<\r_\n$ and $\q$
and $\r$ agree beyond $\n$.  Let $\s\in\U_\I$ have completely
arbitrary entries up to $\s_\m$ and agree with $\q$ beyond $\m$. Since
$\m<\n$, any such $\s$ agrees with $\r$ beyond $\n$ and
$\s_\n=\q_\n<\r_\n$, so $\<\s,\r\>\in\ll\n$.  Since $\q$ and $\s$
agree beyond $\m$, their relationship depends on the relation between
$\q_\m$ and $\s_\m$.  If $\q_\m>\s_\m$ then $\<\q,\s\>\in\cl\m$, and
if $\q_\m<\s_\m$ then $\<\q,\s\>\in\ll\m$.  Both kinds of $\s$ exist,
so $\ll\n\subseteq\cl\m|\ll\n\cap\ll\m|\ll\n$.  From the two
inclusions we have proved, it follows that
$\ll\n=\cl\m|\ll\n=\ll\m|\ll\n$. The other equations in \eqref{rule1}
and \eqref{rule2} can be proved similarly.

Rules \eqref{rule-3}, \eqref{rule-2}, \eqref{rule-1}, \eqref{rule0},
\eqref{rule1}, and \eqref{rule2} show that relative products of atoms
are unions of atoms, hence the set of unions of sets of atoms is
closed under relative multiplication, completing the proof of Theorem
\ref{lemma1}.
\endproof
\begin{lemma}\label{commut}
  Relative multiplication is commutative in $\sra_\I$.
\end{lemma}\proof
Relative multiplication distributes over arbitrary unions, so if
$\mathcal\X,\mathcal\Y\subseteq\At_\I$ then
\begin{align*}
  \bigcup\mathcal\X|\bigcup\mathcal\Y
  &=\bigcup\{\X|\Y\colon\X\in\mathcal\X,\,\Y\in\mathcal\Y\}\\
  &=\bigcup\{\Y|\X\colon\X\in\mathcal\X,\,\Y\in\mathcal\Y\}
  &&\text{\eqref{rule-2}, \eqref{rule0}--\eqref{rule2}}\\
  &=\bigcup\mathcal\Y|\bigcup\mathcal\X.
\end{align*}\endproof
Chains constructed in the next section will be shown in Theorem
\ref{iso} to be Sugihara chains by means of the following
computational rules.
\begin{lemma} Let $\I\subseteq\Zn$. For all $\n,\m\in\I$,
  \begin{align}
    \label{s4a}
    \ll{(-\infty,\n]}|\ll{(-\infty,\m]}
    &=\ll{(-\infty,\n]}\cup\ll{(-\infty,\m]},\\
    \label{s6a}
    \cl{(-\infty,\n]}|\cl{(-\infty,\m]}
    &=\cl{(-\infty,\n]}\cup\cl{(-\infty,\m]},\\
    \label{s5}
    \ll{[\n,\infty)}|\ll{[\m,\infty)}
    &=\ll{[\n,\infty)}\cap\ll{[\m,\infty)},\\
    \label{s7}
    \cl{[\n,\infty)}|\cl{[\m,\infty)}
    &=\cl{[\n,\infty)}\cap\cl{[\m,\infty)},\\
    \label{s7a}
    \cl{(-\infty,\infty)}|\cl{[\m,\infty)}&=\cl{[\m,\infty)},\\
    \label{s7b}
    \cl{(-\infty,\infty)}|\ll{(-\infty,\m]}
    &=\cl{(-\infty,\infty)}\cup\Id_\I\cup\ll{(-\infty,\m]},\\
    \label{s9a}
    \text{if $\n<\m$ then}\quad\ll{(-\infty,\n]}|\cl{[\m,\infty)}
    &=\cl{[\m,\infty)},\\
    \label{s10b}\text{if $\n\geq\m$ then}\quad\ll{(-\infty,\n]}|\cl\m
    &=\cl{(-\infty,\m]}\cup\Id_\I\cup\ll{(-\infty,\n]},\\
    \label{s8a}
    \ll{(-\infty,\n]}|\cl{[\m,\infty)} &=\begin{cases}
      \cl{[\m,\infty)}&\text{ if $\n<\m$,}\\
      \cl{(-\infty,\infty)}\cup\Id_\I\cup\ll{(-\infty,\n]} &\text{ if
        $\n\geq\m$.}
    \end{cases}
  \end{align}
\end{lemma}\proof
In the computations proving \eqref{s4a}--\eqref{s8a} we use
\eqref{rule-3}, \eqref{rule-2}--\eqref{rule2}, and the fact that
relative multiplication distributes over arbitrary unions of
relations.  \eqref{s4a} holds because
\begin{align*}
  \ll{(-\infty,\n]}|\ll{(-\infty,\m]} &=\bigcup\{\ll\k|\ll\ell
  \colon{\n\geq\k\in\I,\,\m\geq\ell\in\I}\}\\
  &=\bigcup\{\ll{\Max(\k,\ell)}\colon{\n\geq\k\in\I,
    \,\m\geq\ell\in\I}\}
  &&\text{\eqref{rule1}}\\
  &=\ll{(-\infty,\Max(\n,\m)]}\\
  &=\ll{(-\infty,\n]}\cup\ll{(-\infty,\m]}.
\end{align*}
Taking converses of both sides in \eqref{s4a} gives \eqref{s6a}. For
\eqref{s5},
\begin{align*}
  \ll{[\n,\infty)}|\ll{[\m,\infty)}
  &=\bigcup\{\ll\k|\ll\ell\colon
  {\n\leq\k\in\I,\,\m\leq\ell\in\I}\}\\
  &=\bigcup\{\ll{\Max(\k,\ell)}\colon
  {\n\leq\k\in\I,\,\m\leq\ell\in\I}\}
  &&\text{\eqref{rule1}}\\
  &=\ll{[\Max(\n,\m),\infty)}\\
  &=\ll{[\n,\infty)}\cap\ll{[\m,\infty)}.
\end{align*}
Applying conversion to \eqref{s5} gives \eqref{s7}. For \eqref{s7a},
\begin{align*}
  \cl{(-\infty,\infty)}|\cl{[\m,\infty)}
  &=\bigcup\{\cl\k|\cl\ell\colon
  {\k\in\I,\,\m\leq\ell\in\I}\}\\
  &=\bigcup\{\cl{\Max(\k,\ell)}\colon
  {\k\in\I,\,\m\leq\ell\in\I}\}
  &&\text{\eqref{rule2}}\\
  &=\cl{[\m,\infty)}.
\end{align*}
For \eqref{s7b}, by \eqref{notation1}, \eqref{notation2}, and
distributivity we have
\begin{align}\label{have}
  \cl{(-\infty,\infty)}|\ll{(-\infty,\m]}=\bigcup\{\cl\k|\ll\ell
  \colon{\k\in\I,\,\m\geq\ell\in\I}\}.
\end{align}
Assume $\k\in\I$ and $\m\geq\ell\in\I$.  If
$\k\neq\ell$ then 
\begin{align*}
  \cl\k|\ll\ell&\subseteq\cl\k\cup\ll\ell
  &&\text{\eqref{rule1}, \eqref{rule2}}\\
  &\subseteq\cl{(-\infty,\infty)}\cup\Id_\I\cup\ll{(-\infty,\m]},
\end{align*}
while if $\k=\ell\leq\m$ then
\begin{align*}
  \cl\k|\ll\ell%
  &=\cl\k|\ll\k\\
  &=\cl{(-\infty,\k]}\cup\Id_\I\cup\ll{(-\infty,\k]}
  &&\text{\eqref{rule0}}\\
  &\subseteq\cl{(-\infty,\infty)}\cup\Id_\I\cup\ll{(-\infty,\m]}.
\end{align*}
Along with \eqref{have}, this shows
\begin{align}\label{goal}
  \cl{(-\infty,\infty)}|\ll{(-\infty,\m]}
  &\subseteq\cl{(-\infty,\infty)}\cup\Id_\I\cup\ll{(-\infty,\m]}.
\end{align}
For the other direction, note that
$\cl\m\subseteq\cl{(-\infty,\infty)}$ and
$\ll\m\subseteq\ll{(-\infty,\m]}$ since $\m\in\I$, hence, by
\eqref{rule0},
\begin{align*}
  \cl{(-\infty,\infty)}|\ll{(-\infty,\m]} \supseteq \cl\m|\ll\m%
  =\cl{(-\infty,\m]}\cup\Id_\I\cup\ll{(-\infty,\m]}.
\end{align*}
What remains is to show
$\cl{(-\infty,\infty)}|\ll{(-\infty,\m]}\supseteq\cl\k$ whenever
$\k\in\I$ and $\k>\m$.  From $\k,\m\in\I$ we get
$\cl\k\subseteq\cl{(-\infty,\infty)}$ and
$\ll\m\subseteq\ll{(-\infty,\m]}$, so
\begin{align*}
  \cl\k&=\cl\k|\ll\m&&\text{$\k>\m$}\\
  &\subseteq\cl{(-\infty,\infty)}|\ll{(-\infty,\m]},
\end{align*}
completing the proof of \eqref{s7b}.  For \eqref{s9a}, if $\n<\m$ then
\begin{align*}
  \ll{(-\infty,\n]}|\cl{[\m,\infty)}
  &=\bigcup\{\ll\k|\cl\ell\colon
  {\k\in\I,\,\k\leq\n<\m\leq\ell\in\I}\}\\
  &=\bigcup\{\cl\ell\colon{\k\in\I,\,\k\leq\n<\m\leq\ell\in\I}\}
  &&\eqref{rule2}\\
  & =\cl{[\m,\infty)}.
\end{align*}
For \eqref{s10b}, if $\n\geq\m$ then, by \eqref{rule0}--\eqref{rule2},
\begin{align*}
  &\ll{(-\infty,\n]}|\cl\m=\bigcup_{\n\geq\k\in\I}\ll\k|\cl\m\\
  &= \bigg(\bigcup_{\n\geq\k\in\I,\,\k<\m}\ll\k|\cl\m\bigg) \cup
  \big(\ll\m|\cl\m\big) \cup
  \bigg(\bigcup_{\n\geq\k\in\I,\,\k>\m}\ll\k|\cl\m\bigg)
  \\
  &= \bigg(\bigcup_{\n\geq\k\in\I,\,\k<\m}\cl\m\bigg) \cup
  \big(\cl{(-\infty,\m]}\cup\Id_\I\cup\ll{(-\infty,\m]}\big) \cup
  \bigg(\bigcup_{\n\geq\k\in\I,\,\k>\m}\ll\k\bigg)
  \\
  &=\cl{(-\infty,\m]}\cup\Id_\I\cup\ll{(-\infty,\n]}.
\end{align*}
Finally we prove \eqref{s8a}.  The first case, in which $\n<\m$, holds
by \eqref{s9a}.  If $\n\geq\m$ then
\begin{align*}
  \ll{(-\infty,\n]}|\cl{[\m,\infty)}
  &=\ll{(-\infty,\n]}|\cl{[\m,\n]}\cup\ll{(-\infty,\n]}|\cl{[\n+1,\infty)}\\
  &=\ll{(-\infty,\n]}|\cl{[\m,\n]}\cup\cl{[\n+1,\infty)}
  &&\text{\eqref{s9a}}\\
  &=\bigg(\bigcup_{\m\leq\ell\leq\n,\,\ell\in\I}\ll{(-\infty,\n]}|\cl\ell\bigg)
  \cup \cl{[\n+1,\infty)}\\
  &=\bigcup_{\m\leq\ell\leq\n,\,\ell\in\I}\left(\Id_\I
    \cup\ll{(-\infty,\n]}\cup \cl{(-\infty,\ell]}\right)
  \cup\cl{[\n+1,\infty)}&&\text{\eqref{s10b}}\\
  &=\ll{(-\infty,\n]}\cup\Id_\I\cup\cl{(-\infty,\n]}\cup\cl{[\n+1,\infty)}\\
  &=\cl{(-\infty,\infty)}\cup\Id_\I\cup\ll{(-\infty,\n]}.
\end{align*}\endproof
\section{Sugihara chains}\label{sect6}
\begin{definition}\label{chain}
  For every $\I\subseteq\Zn$, let
  \begin{align*}
    \Ch_\I&=\{\SI\n\colon-\n\in\I\}\cup\{\TI\n\colon\n\in\I\},\\
    \Ch'_\I&=\{\SI\n\colon-\n\in\I\}\cup\{\TH\n\colon\n\in\I\},
  \end{align*}
  where, for every $\n\in\Zn$,
  \begin{align*}
    \SI\n&=\cl{[-\n,\infty)},&
    \TI\n&=\cl{(-\infty,\infty)}\cup\Id_\I\cup\ll{(-\infty,\n-1]},\\
    &&\TH\n&=\cl{(-\infty,\infty)}\cup\ll{(-\infty,\n-1]}.
  \end{align*}
\end{definition}
It follows from \eqref{notation1} and \eqref{notation2} that the
relations in $\Ch_\I$ and $\Ch'_\I$ form chains under inclusion.  They
are shown in Theorem \ref{iso} to be the universes of Sugihara chains.
When $\I=\Zn$, the order types of $\Ch_\Zn$ and $\Ch'_\Zn$ are the
same as Sugihara's original $\Su_{\Zn+\Zn}$, and the resulting
Sugihara chains are both isomorphic to $\Su_{\Zn+\Zn}$:
\begin{align*}
\tag{$\Ch_\Zn$}\cdots&\subseteq%
  \SI{-2}	\subseteq%
  \SI{-1}	\subseteq%
  \SI0		\subseteq%
  \SI1		\subseteq%
  \cdots&&
  \cdots 	\subseteq%
  \TI{-1}	\subseteq%
  \TI0		\subseteq%
  \TI1		\subseteq%
  \TI2		\subseteq%
  \cdots,\\
\tag{$\Ch'_\Zn$}\cdots&\subseteq%
  \SI{-2}	\subseteq%
  \SI{-1}	\subseteq%
  \SI0		\subseteq%
  \SI1		\subseteq%
  \cdots&&
  \cdots 	\subseteq%
  \TH{-1}	\subseteq%
  \TH0		\subseteq%
  \TH1		\subseteq%
  \TH2		\subseteq%
  \cdots.
\end{align*}
In $\Su_{\Zn+\Zn}$, the designated elements are the ones in the
second, larger copy of $\Zn$.  In the Sugihara chain with universe
$\Ch_\Zn$, the designated relations are the ones that contain the
identity relation on $\U_\Zn$, but all the relations in $\Ch'_\Zn$ are
disjoint from the identity relation.

If $\I=\emptyset$ then $\U_\emptyset$ is a singleton containing just
the function that is constantly zero, and $\Ch_\emptyset =
\Ch'_\emptyset = \emptyset$.  If $\I=\{0\}$ then $\sra_{\{0\}}$ is
isomorphic to Belnap's relation algebra and
\begin{align*}
  \Ch_{\{0\}}&=\{\S^{\{0\}}_0,\T^{\{0\}}_0\},\qquad
  &\Ch'_{\{0\}}&=\{\S^{\{0\}}_0\},\\
  \S^{\{0\}}_0&=\cl0={\hat\T}^{\{0\}}_0,\qquad
  &\T^{\{0\}}_0&=\cl0\cup\Id_{\{0\}}.
\end{align*} 
We thus obtain the two Sugihara chains $\Ch_{\{0\}} = \{\cl0,\cl0 \cup
\Id_{\{0\}}\}$ and $\Ch'_{\{0\}} = \{\cl0\}$, which match up with the
Sugihara chains $\{<,\leq\}$ and $\{<\}$ in Table \ref{ten}, under the
isomorphism $\f$ defined after Definition \ref{defU}.  For a final
example, if $\I=\{0,1\}$, then
\begin{align*}
  \Ch_{\{0,1\}}&=\{\S^{\{0,1\}}_{-1},\S^{\{0,1\}}_0,
  \T^{\{0,1\}}_0,\T^{\{0,1\}}_1\},
  &\Ch'_{\{0,1\}}&=\{\S^{\{0,1\}}_{-1},{\hat\T}^{\{0,1\}}_0,
  {\hat\T}^{\{0,1\}}_1\},\\
  \S^{\{0,1\}}_{-1}&=\cl1,&
  \S^{\{0,1\}}_0&={\hat\T}^{\{0,1\}}_0=\cl0\cup\cl1,\\
  \T^{\{0,1\}}_0&=\cl0\cup\cl1\cup\Id_{\{0,1\}},&
  \T^{\{0,1\}}_1 &=\ll0\cup\cl0\cup\cl1\cup\Id_{\{0,1\}},\\
  &&{\hat\T}^{\{0,1\}}_1 &=\ll0\cup\cl0\cup\cl1.
\end{align*}
Note that $\Ch_{\{0,1\}}$ and $\Ch'_{\{0,1\}}$ can be extended by
adding the empty relation at one end, and the universal relation to
$\Ch_{\{0,1\}}$, or the diversity relation to $\Ch'_{\{0,1\}}$, at the
other end (or both), thus creating Sugihara chains of sizes 5 and 6.
There are four relations in $\sra_{\{0,1\}}$ that are fixed by
$\rmin'$, namely $\ll0\cup\ll1$, $\cl0\cup\ll1$, $\ll0\cup\cl1$, and
$\cl0\cup\cl1$.  Two of these relations appear in the middle of two
Sugihara chains of length 5. The union of these two chains forms a
definitional reduct of $\sra_{\{0,1\}}$ that is isomorphic to the
crystal lattice in \SS\ref{sect8}.
\begin{theorem}\label{iso}
  For every $\I\subseteq\Zn$, $\<\Ch_\I,\cup,\cap,\to,\rmin\>$ and
  $\<\Ch'_\I,\cup,\cap,\to',\rmin'\>$ are Sugihara chains.  In
  particular, $\<\Ch_\Zn,\cup,\cap,\to,\rmin\>$ and
  $\<\Ch'_\Zn,\cup,\cap,\to',\rmin'\>$ are isomorphic to the original
  $\Su_{\Zn+\Zn}$.
\end{theorem}\proof 
By \eqref{rule-3} and \eqref{notation2},
\begin{align*}
  \rmin\SI\n&=\min{\conv{\SI\n}} =\min{\conv{\cl{[-\n,\infty)}}}
  =\min{\ll{[-\n,\infty)}} =\cl{(-\infty,\infty)}\cup\Id_\I\cup
  \ll{(-\infty,-\n-1]}=\TI{-\n},
\end{align*}
so
\begin{align}\label{rmin'rule}
  \rmin'\SI\n&=\rmin\SI\n\cap\Di_\I=\TI{-\n}\cap\Di_\I=\TH{-\n}.
\end{align}
It is straightforward to verify that $\rmin$ is an order-reversing
involution on all relations, and that $\rmin'$ is an order-reversing
involution on relations included in $\Di_\I$.  Therefore we have
$\rmin\TI{-\n}=\SI\n$, and $\rmin'(\TH{-\n})=\SI\n$ since
$\SI\n\cup\TH{-\n}\subseteq\Di_\I$.  It follows that $\Ch_\I$ and
$\Ch'_\I$ are closed under converse-complementation $\rmin$ and
relativized converse-complementation $\rmin'$, respectively.  Since
$\rmin\emptyset=(\U_\I)^2$ and $\rmin'\emptyset=\Di_\I$,
converse-complementation and relativized converse-complementation are
also order-reversing involutions on $\Ch_\I\cup \{\emptyset,
(\U_\I)^2\}$ and $\Ch'_\I\cup\{\emptyset,\Di_\I\}$, respectively.
Turning to relative products, we show for all $\n,\m\in\I$,
\begin{align}
  \label{p1}
  \SI\n|\SI\m&=\SI\n\cap\SI\m,\\
  \label{p2}
  \SI\n|\TI\m=\TI\m|\SI\n&=\begin{cases}\SI\n&\text{ if $\n\leq-\m$,}\\
    \TI\m &\text{ if $\n>-\m$,}\end{cases}\\
  \label{p3} \TI\m|\TI\n&=\TI\m\cup\TI\n,  \\
  \label{pq2}
  \SI\n|\TH\m=\TH\m|\SI\n&=\begin{cases}\SI\n
    &\text{ if $\n\leq-\m$,}\\
    \TH\m\cup\Id_\I &\text{ if $\n>-\m$,}\end{cases}\\
  \TH\m|\TH\n&=\TH\m\cup\TH\n\cup\Id_\I.
  \label{pq3}
\end{align}
For \eqref{p1} and \eqref{p2} we have
\begin{align*}
  \SI\n|\SI\m&=\cl{[-\n,\infty)}|\cl{[-\m,\infty)}\\
  &=\cl{[-\n,\infty)}\cap\cl{[-\m,\infty)}&&\eqref{s7}\\
  &=\SI\n\cap\SI\m,\\
  \TI\m|\SI\n &=\left(\cl{(-\infty,\infty)}\cup\Id_\I\cup
    \ll{(-\infty,\m-1]}\right)|\cl{[-\n,\infty)}\\
  &=\cl{(-\infty,\infty)}|\cl{[-\n,\infty)}\cup
  \Id_\I|\cl{[-\n,\infty)}\cup
  \ll{(-\infty,\m-1]}|\cl{[-\n,\infty)}\\
  &=\cl{[-\n,\infty)}\cup\cl{[-\n,\infty)}
  \cup\ll{(-\infty,\m-1]}|\cl{[-\n,\infty)}
  &&\text{\eqref{s7a}, \eqref{rule-2}}\\
  &=\begin{cases}\cl{[-\n,\infty)}&\text{ if $\m-1<-\n$}\\
    \cl{(-\infty,\infty)}\cup\Id_\I\cup\ll{(-\infty,\m-1]} &\text{
      if $\m-1\geq-\n$}\end{cases} &&\text{\eqref{s8a}}\\
  &=\begin{cases}\SI\n&\text{ if $\n\leq-\m$,}\\
    \TI\m &\text{ if $\n>-\m$.}\end{cases}
\end{align*}
For \eqref{p3} we start with the observation that
\begin{align*}
  \TI\m|\TI\n&=
  \left(\cl{(-\infty,\infty)}\cup\Id_\I\cup\ll{(-\infty,\m-1]}\right)|
  \left(\cl{(-\infty,\infty)}\cup\Id_\I\cup\ll{(-\infty,\n-1]}\right).
\end{align*} 
Multiplying this out yields these nine relative products.
\begin{align*}
  \cl{(-\infty,\infty)}|\cl{(-\infty,\infty)} &=\cl{(-\infty,\infty)},
  &&\text{\eqref{rule-1}}\\
  \cl{(-\infty,\infty)}|\Id_\I &=\cl{(-\infty,\infty)},
  &&\text{\eqref{rule-2}}\\
  \cl{(-\infty,\infty)}|\ll{(-\infty,\n-1]},
  &=\cl{(-\infty,\infty)}\cup\Id_\I\cup\ll{(-\infty,\n-1]},
  &&\text{\eqref{s7b}}\\
  \Id_\I|\cl{(-\infty,\infty)} &=\cl{(-\infty,\infty)},
  &&\text{\eqref{rule-2}}\\
  \Id_\I|\Id_\I &=\Id_\I,
  &&\text{\eqref{rule-2}}\\
  \Id_\I|\ll{(-\infty,\n-1]} &=\ll{(-\infty,\n-1]},
  &&\text{\eqref{rule-2}}\\
  \ll{(-\infty,\m-1]}|\cl{(-\infty,\infty)}
  &=\cl{(-\infty,\infty)}\cup\Id_\I\cup\ll{(-\infty,\m-1]},
  &&\text{\eqref{s7b}, Lemma \ref{commut}}\\
  \ll{(-\infty,\m-1]}|\Id_\I &=\ll{(-\infty,\m-1]},
  &&\text{\eqref{rule-2}}\\
  \ll{(-\infty,\m-1]}|\ll{(-\infty,\n-1]}
  &=\ll{(-\infty,\m-1]}\cup\ll{(-\infty,\n-1]}. &&\text{\eqref{s4a}}
\end{align*}
Taking the union of the relations on the right gives us
\begin{align*}
  \TI\m|\TI\n&= \cl{(-\infty,\infty)}\cup\Id_\I\cup
  \ll{(-\infty,\m-1]}\cup\ll{(-\infty,\n-1]}=\TI\m\cup\TI\n,
\end{align*}
so \eqref{p3} holds.  The proofs of \eqref{pq2} and \eqref{pq3} are
somewhat simpler. They can be obtained from the computations just
given by deleting references to $\Id_\I$ on the left sides of the
equations, and expressing the relations on the right sides in terms of
$\Id_\I$ and the relations in $\Ch'_\I$.  Recall from Definition
\ref{def1} that in a Sugihara chain, $\to$ is defined by
\begin{align*}
  \A\to\B &=\begin{cases}
    \rmin\A\lor\B&\text{if $\A\leq\B$,}\\
    \rmin\A\land\B&\text{if $\A>\B$.}\end{cases}
\end{align*}
Substitute $\rmin\B$ for $\B$ and apply $\rmin$ to both sides. The
double negation and De Morgan laws for $\land$, $\lor$, and $\rmin$
hold in every Sugihara chain, so
\begin{align*} \rmin(\A\to\rmin\B) &=\begin{cases}
    \A\land\B&\text{if $\A\leq\rmin\B$,}\\
    \A\lor\B&\text{if $\A>\rmin\B$.}\end{cases}
\end{align*}
To show residuation and relativized residuation act like Sugihara's
$\to$, we will use the latter equation.  By some elementary
calculations starting from the definitions in Table \ref{table1} of
relative multiplication, converse-complementation, residuation and
their relativized counterparts, we get
\begin{align}\label{fuse}
  \rmin(\A\to\rmin\B)=\B|\A,\qquad\qquad
  \rmin'(\A\to'\rmin'\B)=\B|'\A,
\end{align} 
hence all we need to show is that for any $\A,\B\in\Ch_\I$,
\begin{align}
  \label{show}
  \B|\A&=\begin{cases}
    \A\cap\B&\text{if $\A\subseteq\rmin\B$,}\\
    \A\cup\B&\text{if $\A\supset\rmin\B$,}
  \end{cases}
  \intertext{and for any $\A,\B\in\Ch'_\I$,}
  \label{show1}
  \B|'\A&=\begin{cases}
    \A\cap\B&\text{if $\A\subseteq\rmin'\B$,}\\
    \A\cup\B&\text{if $\A\supset\rmin'\B$.}
  \end{cases}
\end{align}
Proof of \eqref{show}. Because of commutativity (Lemma \ref{commut}),
there are just three cases that arise by substituting into
\eqref{show} when $\n,\m\in\I$, $\A\in\{\SI\n,\TI\n\}$, and
$\B\in\{\SI\m,\TI\m\}$.

{\bf Case~1}. $\A=\SI\n$, $\B=\SI\m$, $\rmin\B=\TI{-\m}$.  The first
case in \eqref{show} applies because $\SI\n\subseteq\TI{-\m}$.  By
\eqref{p1}, $\B|\A=\SI\m|\SI\n=\SI\n\cap\SI\m=\A\cap\B$. This agrees
with \eqref{show}, and shows that \eqref{show} holds.

{\bf Case~2}. $\A=\TI\n$, $\B=\TI\m$, $\rmin\B=\SI{-\m}$.  The second
case in \eqref{show} applies since $\TI\n\supset\SI{-\m}$.  By
\eqref{p3} we have $\B|\A=\TI\m|\TI\n=\TI\n\cup\TI\m=\A\cup\B$, as
required for \eqref{show} to hold.

{\bf Case~3}. $\A=\SI\n$, $\B=\TI\m$, $\rmin\B=\SI{-\m}$. Since
$\A\subset\B$ in this case, \eqref{show} simplifies into the form
proved below.
\begin{align*}
  \B|\A=\TI\m|\SI\n &=\begin{cases}
    \SI\n&\text{if $\n\leq-\m$}\\
    \TI\m&\text{if $\n>-\m$}
  \end{cases}&&\text{\eqref{p2}}\\
  &=\begin{cases}
    \SI\n&\text{if $\SI\n\subseteq\SI{-\m}$}\\
    \TI\m&\text{if $\SI\n\supset\SI{-\m}$}
  \end{cases}&&\text{}\\
  &=\begin{cases}
    \A&\text{if $\A\subseteq\rmin\B$,}\\
    \B&\text{if $\A\supset\rmin\B$.}
  \end{cases}
\end{align*}
Proof of \eqref{show1}. Again there are three cases.

{\bf Case~1}. $\A=\SI\n$, $\B=\SI\m$, $\rmin'\B=\TH{-\m}$.  The first
case in \eqref{show1} applies because $\SI\n\subseteq\TH{-\m}$.  By
\eqref{p1} and $\SI\n\cup\SI\m\subseteq\Di_\I$,
$\B|'\A=\SI\m|\SI\n\cap\Di_\I=\SI\n\cap\SI\m=\A\cap\B$, which agrees
with \eqref{show1}.

{\bf Case~2}. $\A=\TH\n$, $\B=\TH\m$, $\rmin'\B=\SI{-\m}$. By
$\TH\n\supset\SI{-\m}$, the second case in \eqref{show1} applies.  By
\eqref{p3} and $\TH\n\cup\TH\m\subseteq\Di_\I$,
$\B|'\A=\TH\m|\TH\n\cap\Di_\I=\TH\n\cup\TH\m=\A\cup\B$, so
\eqref{show1} holds.

{\bf Case~3}. $\A=\SI\n$, $\B=\TH\m$, $\rmin'\B=\SI{-\m}$. In this
case $\A\subset\B$, so for \eqref{show1} we need only show
\begin{align*}
  \B|'\A=\TH\m|\SI\n\cap\Di_I&=\begin{cases}
    \SI\n&\text{if $\n\leq-\m$}\\
    \TH\m&\text{if $\n>-\m$}
  \end{cases}&&\text{\eqref{pq2}, $\TH\m\cup\SI\n\subseteq\Di_\I$}\\
  &=\begin{cases}
    \A&\text{if $\A=\SI\n\subseteq\SI{-\m}=\rmin'\B$,}\\
    \B&\text{if $\A=\SI\n\supset\SI{-\m}=\rmin'\B$.}
  \end{cases}
\end{align*}
This completes the proof that $\<\Ch_\I,\cup,\cap,\to,\rmin\>$ and
$\<\Ch'_\I,\cup,\cap,\to',\rmin'\>$ are Sugihara chains. Because of
the match between the order types when $\I=\Zn$, illustrated in the
remarks preceding Theorem \ref{iso}, we also conclude that
\begin{align*}
  \Su_{\Zn+\Zn} \cong \<\Ch_\Zn,\cup,\cap,\to,\rmin\> \cong
  \<\Ch'_\Zn,\cup,\cap,\to',\rmin'\>,
\end{align*}
which completes the proof of Theorem \ref{iso}.
\endproof
The set of converses of a Sugihara chain is another Sugihara chain.
Applying this observation to $\Ch_\I$, we let
\begin{align*}
  {\con\S}^\I_\n&=\ll{[-\n,\infty)},\\
  {\con\T}^\I_\n&=\ll{(-\infty,\infty)}\cup\Id_\I\cup\cl{(-\infty,\n-1]},\\
  {\con\Ch}_\I&=\{{\con\S}^\I_\n\colon-\n\in\I\}\cup
  \{{\con\T}^\I_\n\colon\n\in\I\}.
\end{align*}
Then ${\con\Ch}_\I$ is the other copy of the Sugihara chain $\Ch_\I$
in $\sra_\I$.  Observe that
\begin{align*}
  \mathcal{X}_\I=\Ch_\I\cup{\con\Ch}_\I\cup
  \{\Id_\I,\Di_\I,\emptyset,\U^2\}
\end{align*}
is closed under union, intersection, converse-complementation, and
residuation. In addition, $\mathcal{X}_\I$ is closed and commutative
under relative multiplication. All the relations in $\mathcal{X}_\I$
are dense.  The only non-transitive relation in $\mathcal{X}_\I$ is
$\Di_\I$.  When $\I=\{0\}$, $\mathcal\X_{\{0\}}$ coincides with the
entire universe of $\sra_{\{0\}}$, reflecting the fact that the
diversity relation $\neq$ is the only non-transitive relation in
$\M_0$.
\section{Reducts, relation algebras, and %
  atom structures}\label{sect7}
\begin{definition}\label{reducts}
  A {\bf definitional reduct} of an algebra $\mathfrak\A$ is obtained
  by omitting some of the fundamental operations of $\mathfrak\A$ and
  adding some operations that are term-definable in $\mathfrak\A$.  A
  {\bf definitional subreduct} is a subalgebra of a definitional
  reduct.
\end{definition}
For example, when defined as in \eqref{belnap}, Belnap's $\MM$ is a
definitional reduct of Belnap's relation algebra $\gc\M_0$, but not
conversely.  They both have the same universe, but $\gc\M_0$ has
operations not definable from the operations of $\MM$.  We use two
methods to obtain definitional subreducts, called {\bf direct} and
{\bf relativized}.  They apply to all relation algebras, although we
will be primarily interested in applying them to proper relation
algebras (Definition \ref{proper}), so we review basic definitions and
facts about relation algebras.  Good resources for relation algebras
are \cite{MR3699802, MR3699801, MR1935083, MR2269199}, especially the
first two.
\begin{definition}\label{radef}
  A {\bf relation algebra} is an algebra $\gc\A=\<\A,+,\cdot,
  \min\blank, 0, 1, \rp, \con\blank, \id\>$, consisting of a set $\A$,
  binary operations $+$ and $\rp$ on $\A$, unary operations
  $\min\blank$ and $\con\blank$ on $\A$, and a distinguished element
  $\id\in\A$, called the {\bf identity element} of $\gc\A$, such that
  $\<\A,+,\cdot,\min\blank,0,1\>$ is a Boolean algebra and $\gc\A$
  satisfies the axioms %
\begin{align*}
  \tag{r1} (\x\rp\y)\rp\z&=\x\rp(\y\rp\z),\\
  \tag{r2} (\x+\y)\rp\z&=(\x\rp\z)+(\x\rp\z),\\
  \tag{r3} \x&=\x\rp\id,\\
  \tag{r4} \con{\con\x}&=\x,\\
  \tag{r5} \con{\x+\y}&=\con\x+\con\y,\\
  \tag{r6} \con{\x\rp\y}&=\con\y\rp\con\x,\\
  \tag{r7} \con\x\rp\min{\x\rp\y}&\leq\min\y,
\end{align*}
where $\x\leq\y$ iff $\x+\y=\y$.  Define the {\bf diversity element}
by $\di=\min\id$. An element $\a\in\A$ is an {\bf atom} if $\a\neq0$
and for all $\x\in\A$, if $\x\leq\a$ then $\x=0$ or $\x=\a$.  $\gc\A$
is {\bf atomic} if for every non-zero element $\x\in\A$ there is an
atom $\a\in\A$ such that $\a\leq\x$. An element $\x\in\A$ is {\bf
  symmetric} if $\con\x=\x$, {\bf dense} if $\x\leq\x\rp\x$, and {\bf
  transitive} if $\x\rp\x\leq\x$.  The algebra $\gc\A$ is {\bf
  symmetric} if all its elements are symmetric, {\bf dense} if its
elements are dense, {\bf commutative} if it satisfies $\x\rp\y =
\y\rp\x$, and {\bf Boolean} if $\id=1$.
\end{definition}
Proper relation algebras are, indeed, relation algebras.  Boolean
relation algebras are symmetric and also satisfy $\x\rp\y=\x\cdot\y$.
Their name is based on the observation, made after \cite[Theorem
4.35]{MR44502}, that if $\<\A,+,\cdot,\min\blank,0,1\>$ is a Boolean
algebra and $\con\x$ is defined to be $\x$, then $\<\A,+, \cdot,
\min\blank, 0, 1, \cdot, \con\blank,1\>$ is a relation algebra.  Each
of the identities $\id=1$ and $\x\rp\y=\x\cdot\y$ characterizes
Boolean relation algebras \cite[Lemma 3.1]{MR3699801}.  Boolean
relation algebras are representable \cite[Theorem 17.5]{MR3699802}.  A
relation algebra $\gc\A$ is {\bf simple} (has no non-trivial
homomorphic images) if and only if $1\rp\x\rp1=1$ whenever
$0\neq\x\in\A$ \cite[Theorem 4.10]{MR44502}.  A relation algebra is
{\bf integral} (has no zero divisors) if and only if $\id$ is an atom
\cite[Theorem 4.17]{MR44502}. In every relation algebra, the {\bf
  converse} $\con\a$ of an atom $\a$ is an atom \cite[Theorem
4.3(xii)]{MR44502}, which allows the following definition of atom
structure. The definition of complex algebra is the special case of
\cite[Definition 3.8]{MR44502} that applies to relation algebras.
\begin{definition}{\cite[Definitions 2.1, 3.2]{MR662049}}
  \label{cmat}
  \begin{enumerate}
  \item If $\gc\U=\<\U,\R,\f,\I\>$ is a structure where
    $\R\subseteq\U^3$, $\f\colon\U\to\U$, and $\I\subseteq\U$, then
    the {\bf complex algebra} of $\gc\U$ is 
    \begin{equation*}
      \Cm{\gc\U}=\<\wp(\U), \cup, \cap, \min\blank, \emptyset, \U, 
      \rp, \con\blank, \I\>,
    \end{equation*}
    where $\wp(\U)$ is the powerset of $\U$, $\<\wp(\U), \cup, \cap,
    \min\blank, \emptyset, \U\>$ is the Boolean algebra of all subsets
    of $\U$, and for all $\X,\Y\subseteq\U$, $\con\X =
    \{\con\x:\x\in\X\}$ and $\X\rp\Y = \{\z: \x\rp\y\geq\z\in\at
    \text{ for some $\x\in\X$ and $\y\in\Y$}\}$.
  \item The {\bf atom structure} of an atomic relation algebra $\gc\A$
    is $\<\at,\R,\con\blank,\I\>$, where $\at$ is the set of atoms of
    $\gc\A$, $\R=\{\<\x,\y,\z\>:\x,\y,\z\in\at,\, \x\rp\y\geq\z\}$,
    and $\I=\{\u:\id\geq\u\in\at\}$.
  \end{enumerate}
\end{definition}
The next theorem is the relation algebraic case of \cite[Theorem
3.9]{MR44502}. The specific conditions were first stated earlier in
\cite[\SS4]{MR0037278} in a slightly different but equivalent form.
\begin{theorem}{\cite[Theorem 2.2]{MR662049}}\label{ra}
  The complex algebra $\Cm{\gc\U}$ of a structure $\gc\U = \<\U, \R,
  \f, \I\>$ is a complete and atomic Boolean algebra with operators,
  and $\Cm{\gc\U}$ is a relation algebra if and only if, for all
  $\x,\y,\z\in\U$,
  \begin{enumerate}
  \item\label{refl1} if $\<\x,\y,\z\>\in\R$ then
    $\<\f\x,\z,\y\>\in\R$,
  \item\label{refl2} if $\<\x,\y,\z\>\in\R$ then
    $\<\z,\f\y,\x\>\in\R$,
  \item\label{identity} $\x=\y$ iff there is some $\u\in\I$ such that
    $\<\x,\u,\y\>\in\R$,
  \item\label{pasch} if $\<\v,\w,\x\>,\<\x,\y,\z\>\in\R$ then for some
    $\u\in\U$, $\<\v,\u,\z\>,\<\w,\y,\u\>\in\R$.
  \end{enumerate}
\end{theorem}
It follows from just \ref{refl1}, \ref{refl2}, and \ref{identity} that
$\f$ is an {\bf involution} on $\U$ ($\f\f\x=\x$ for all $\x\in\U$),
$\f\x=\x$ for all $\x\in\I$, $\R$ is the union of {\bf cycles}
\cite[(1), p.\ 710]{MR0037278}, which are sets of the form
\begin{align}\label{cycle}
  [\x,\y,\z]=\{&\<\x,\y,\z\>, \<\z,\f\y,\x\>, \<\f\z,\x,\f\y\>,\\
  \notag &\<\f\y,\f\x,\f\z\>,\<\y,\f\z,\f\x\>, \<\f\x,\z,\y)\>\},
\end{align}
and finally, the complex algebra satisfies axioms
\thetag{r2}--\thetag{r8}. The associative law \thetag{r1} is the only
axiom that may fail, and \thetag{r1} holds if and only if \ref{pasch}
holds. Every relation algebra has a complete and atomic extension,
called its {\bf perfect extension}, {\bf canonical extension}, or {\bf
  canonical embedding algebra} \cite[Theorem 4.21]{MR45086}, from
which we get the following special case of \cite[Theorem 3.10
(Representation Theorem)]{MR44502} that includes the appropriate
conditions for relation algebras.
\begin{theorem}{\cite[Theorems 3.13, 4.3]{MR662049}}\label{repth}
  A relation algebra is complete and atomic if and only if it is
  isomorphic to the complex algebra of its atom structure. An algebra
  $\gc\A=\<\A,+,\cdot,\min\blank,0,1,\rp,\con\blank,\id\>$ is a
  relation algebra if and only if it is isomorphic to a subalgebra of
  the complex algebra of a structure satisfying conditions
  {\normalshape\ref{refl1}}, {\normalshape\ref{refl2}},
  {\normalshape\ref{identity}}, and {\normalshape\ref{pasch}} in
  Theorem {\normalshape\ref{ra}}.
\end{theorem}
For an arbitrary relation algebra $\gc\A$, this structure may be
constructed directly from $\gc\A$ as follows \cite[Theorem
2.11]{MR2616328}.  An {\bf ultrafilter} is a maximal proper subset
$\X\subseteq\A$ such that $\x\cdot\y\in\X$ whenever $\x,\y\in\X$ and
$\x+\y\in\X$ whenever $\x\in\X$ and $\y\in\A$. Let $\U$ be the set of
ultrafilters of $\gc\A$. Let $\R$ be the set of triples $\<\X,\Y,\Z\>$
of ultrafilters such that $\X\rp\Y\subseteq\Z$, let $\f\colon\U\to\U$
be defined by $\f\X=\{\con\x:\x\in\X\}$, and let $\I$ be the set of
ultrafilters that contain $\id$. Then the desired {\bf canonical atom
  structure} is $\<\U,\R,\f,\I\>$. When $\gc\A$ is complete and
atomic, the canonical atom structure is isomorphic to the atom
structure of $\gc\A$.
\begin{definition}\label{methods}
  Let $\gc\A=\<\A,+,\cdot, \min\blank, 0, 1, \rp, \con\blank, \id\>$
  be a relation algebra. Then $\mathfrak\A_r=\<\A,+,\cdot,\to,\rmin\>$
  is the algebra obtained from $\mathfrak\A$ by deleting $\min\blank$,
  $0$, $1$, $\rp$, $\con\blank$, and $\id$, retaining $+$ and $\cdot$,
  and adding operations $\to$ and $\rmin$, defined by
  \begin{align*}
    \x\to\y&=\min{\con\x\rp\min\y},&\rmin\x&=\min{\con\x}, %
    \intertext{and $\mathfrak\A_r'=\<\A,+,\cdot,\to',\rmin'\>$ is the
      algebra obtained from $\mathfrak\A$ by deleting $\min\blank$,
      $0$, $1$, $\rp$, $\con\blank$, and $\id$, retaining $+$ and
      $\cdot$, and adding operations $\to'$ and $\rmin'$, defined by}
    \x\to'\y&=\min{\con{\x\cdot\di}\rp\min{\y\cdot\di}}\cdot\di,&
    \rmin'\x&=\min{\con{\x\cdot\di}}\cdot\di.
  \end{align*}
  $\mathfrak\A_r$ is the {\bf direct reduct} of $\mathfrak\A$, and
  $\mathfrak\A_r'$ is the {\bf relativized reduct} of $\mathfrak\A$.
  An algebra is a {\bf direct subreduct} of $\mathfrak\A$ if it is a
  subalgebra of the direct reduct of $\mathfrak\A$, and a {\bf
    relativized subreduct} of $\mathfrak\A$ if it is a subalgebra of
  the relativized reduct of $\mathfrak\A$.
\end{definition}
When applied to proper relation algebras, the operations $\to$,
$\rmin$, $\to'$, and $\rmin'$ are the ones (with the same names)
defined in Table \ref{table1}. For example, when defined as in
\eqref{belnap}, Belnap's $\MM$ is the direct reduct of the proper
relation algebra $\gc\M_0$. The next theorem is the main result of
this paper.  The part asserting that every finite Sugihara chain of
even cardinality is isomorphic to a direct subreduct of a proper
relation algebra was already proved in \cite[Theorem~6.2]{MR2641636}.
The two innovations that allow us to extend this result to infinite
Sugihara chains and to finite Sugihara chains of odd cardinality are
sequences that are eventually zero and relativized subreducts.
\begin{theorem}\label{main} For every $\I\subseteq\Zn$, $\sra_\I$ is a
  proper relation algebra such that
  \begin{enumerate}
  \item\label{c1} the Sugihara chain $\<\Ch_\I,\cup,\cap,\to,\rmin\>$
    is a direct subreduct of $\sra_\I$, and
  \item\label{c2} the Sugihara chain $\< \Ch'_\I, \cup, \cap, \to',
    \rmin'\>$ is a relativized subreduct of $\sra_\I$.
  \item\label{c3} $\Su_{\Zn^*}$, Sugihara's original $\Su_{\Zn+\Zn}$,
    and all finite Sugihara chains of even cardinality are isomorphic
    to both direct subreducts and relativized subreducts of the proper
    relation algebra $\sra_\Zn$.
  \item\label{c4} $\Su_\Zn$ and all Sugihara chains of odd cardinality
    are isomorphic to relativized subreducts of the proper relation
    algebra $\sra_{\Zn^+}$, where $\Zn^+=\{\n:0<\n\in\Zn\}$.
  \end{enumerate}
\end{theorem}
\proof Recall that $\mathcal\S_\I$ is the universe of the Sugihara
relation algebra $\sra_\I$.  By Definition \ref{defU},
\eqref{notation1}, \eqref{notation2}, and Definition \ref{chain}, we
have $\Ch_\I\subseteq\mathcal\S_\I$ and
$\Ch'_\I\subseteq\mathcal\S_\I$. Furthermore, $\Ch_\I$ is closed under
$\cup$, $\cap$, $\to$, and $\rmin$, and $\Ch'_\I$ is closed under
$\cup$, $\cap$, $\to'$, and $\rmin'$ since $\<\Ch_\I, \cup, \cap, \to,
\rmin\>$ and $\<\Ch'_\I, \cup, \cap, \to', \rmin'\>$ are Sugihara
chains by Theorem \ref{iso}. It follows by Definition \ref{methods}
that $\<\Ch_\I, \cup, \cap, \to, \rmin\>$ is a direct subreduct of
$\sra_\I$ and $\< \Ch'_\I, \cup, \cap, \to', \rmin'\>$ is a
relativized subreduct of $\sra_\I$.  Thus \ref{c1} and \ref{c2} hold.
By Theorem \ref{iso},
\begin{align*}
  \<\Ch_\Zn,\cup,\cap,\to,\rmin\>\cong
  \<\Ch'_\Zn,\cup,\cap,\to',\rmin'\>\cong\Su_{\Zn+\Zn},
\end{align*}
so $\Su_{\Zn+\Zn}$ is isomorphic to both a direct and a relativized
subreduct of $\sra_\Zn$. All of the subalgebras of $\Su_{\Zn+\Zn}$ are
therefore also isomorphic to direct and relativized subreducts of
$\sra_\Zn$.  This includes $\Su_{\Zn^*}$ and all finite Sugihara
chains of even cardinality, thus proving \ref{c3}.

For the Sugihara chains of odd cardinality we proceed differently.
Suppose $\I$ has a minimum element $\m\in\I\subseteq\Zn$. This means
that $\{\n:\m>\n\in\I\} = \emptyset$. Therefore $\ll{(-\infty,\m-1]} =
\bigcup\emptyset = \emptyset$ by \eqref{notation1} and
$\cl{[\m,\infty)} = \cl{(-\infty,\infty)}$ by \eqref{notation2}, so by
\eqref{rmin'rule},
\begin{align*}
  \SI{-\m}=\cl{[\m,\infty)}
  =\cl{(-\infty,\infty)}\cup\ll{(-\infty,\m-1]}=\TH\m=\rmin'(\SI{-\m}).
\end{align*}
Thus the relation $\SI{-\m}$ is fixed by $\rmin'$. If $\I$ is also
infinite, then $\<\Ch'_\I,\cup,\cap,\to',\rmin'\>$ is isomorphic to
$\Su_\Zn$, by an isomorphism that sends the fixed point
$\SI{-\m}=\TH\m$ to $0$, which is the fixed point of negation in
$\Su_\Zn$.  Taking $\I=\Zn^+$, we see that $\Su_\Zn$ is isomorphic to
a relativized subreduct of $\sra_{\Zn^+}$. Every subalgebra of
$\Su_\Zn$ is also isomorphic to a relativized subreduct of
$\sra_{\Zn^+}$. This includes all finite Sugihara chains of odd
cardinality, so \ref{c4} holds.
\endproof
Theorem \ref{main} shows that every finite Sugihara chain $\Su$ is
isomorphic to a chain of binary relations closed under the relevant
operations.  If the identity relation is included in the relations
occurring in the top half of this chain (this is the direct method),
then there cannot be a relation fixed by $\rmin$ and $\Su$ has even
cardinality.  If $\Su$ is odd it can be represented by the relativized
method purely with diversity relations.  The element in the middle of
$\Su$ is mapped to a relation that is its own relativized
converse-complement.  For example, Belnap's relation algebra
$\sra_{\{0\}}$ has two Sugihara chains of length 3, namely
$\{\emptyset,<,\neq\}$ and $\{\emptyset,>,\neq\}$ (see Figure
\ref{fig1} and Table \ref{ten}). Note how the relations in the middle,
namely $<$ and $>$, are fixed by $\rmin'$.  The Sugihara chains of
length 3 are isomorphic to RM3, described on \cite[p.\ 470]{MR0406756}
and \cite[p.\ 92]{MR728950}.  This is yet another algebra of relevance
logic that can be represented as an algebra of binary relations.
Sugihara chains of even cardinality are subalgebras of those with odd
cardinality.  Thus the normal Sugihara chains are isomorphic to
definitional subreducts of proper relation algebras by both the direct
and relativized methods, while non-normal ones need the relativized
method.

From Theorem \ref{main} we know that $\Su_{\Zn^*}$ (``the Sugihara
matrix'' of Anderson and Belnap) is isomorphic to a direct subreduct
of $\sra_{\Zn^+}$.  In \cite{MR2536403} there is a computation
intended to show that this is not possible.  On \cite[p.\
123]{MR2536403},
\begin{quote}
  ``Figure 5 shows some components of the canonical embedding algebra
  of the Sugihara matrix $\Su_{\Zn^*}$ \dots\ Unfortunately, this
  Boolean algebra is not a relation algebra, let alone a transitive or
  a representable one. To show that (r6) is not true, we give a
  concrete counterexample.''
\end{quote} 
The ensuing computation at the bottom of \cite[p.\,122]{MR2536403}
ends with $\{[i):i\leq-2\}$, but should end with $\{[i):i\geq-2\}$.
When corrected, it confirms an instance of axiom (r6) in Definition
\ref{radef}.  It was reasonable to suspect this equation may not hold,
because it corresponds to a property of atom structures of relation
algebras not shared by the model structures of $\RR$.  That property,
identified and called ``tagging'' by Dunn \cite{MR2067967}, says that
if $\<x,y,z\>\in\R$ then $\<fy,fx,fz\>\in\R$. In Theorem \ref{ra},
either one of \ref{refl1} and \ref{refl2} can be replaced by tagging.
None of these three conditions necessarily holds in a model structure
for $\RR$.  Such structures do satisfy the condition that if
$\<x,y,z\>\in\R$ then $\<fz,x,fy\>\in\R$. This condition, together
with any one of of the three conditions \ref{refl1}, \ref{refl2}, and
tagging, can be used in Theorem \ref{ra} instead of \ref{refl1} and
\ref{refl2}, because any of these combinations are enough to prove
that $\R$ is a union of cycles \eqref{cycle}. To obtain relevant model
structures for logics like $\RR$, one must add the conditions
expressing density, that $\<\x,\x,\x\>\in\R$ for all $\x$, and
commutativity, that if $\<\x,\y,\z\>\in\R$ then $\<\y,\x,\z\>\in\R$.

By Theorem \ref{main}, the canonical embedding algebra of the Sugihara
matrix $\Su_{\Zn^*}$ is, in fact, isomorphic to the complete atomic
proper relation algebra $\sra_{\Zn^+}$, hence also isomorphic to the
complex algebra of the canonical atom structure of $\sra_{\Zn^+}$.
Although $\sra_{\Zn^+}$ is commutative and dense, not all of its
elements are transitive.  For example, the diversity relation
$\Di_{\Zn^+}$ is not transitive.  However, $\sra_{\Zn^+}$ does have
subsets that contain only transitive (and dense) relations and are
closed under the relevant operations. As we have seen, a copy of the
Sugihara matrix $\Su_{\Zn^*}$ is among them.
\section{The crystal lattice, Church's diamond, %
  and Meyer's RM84}\label{sect8}
\subsection*{The crystal lattice}
\begin{figure}
  \begin{picture}(40,32)(0,0)%
    \put(20,30){\circle*{1.1}}%
    \put(21,30){$\Di_{\{0,1\}}=\rmin'(\emptyset)$}%
    \put(20,20){\circle*{1.1}}%
    \put(21.5,20){$\ll0\cup\cl0\cup\ll1=\rmin'(\ll1)$}%
    \put(10,15){\circle*{1}}%
    \put(-8.5,15){$\rmin'(\ll0\cup\ll1)=\ll0\cup\ll1$}%
    \put(30,15){\circle*{1}{$\cl0\cup\ll1=\rmin'(\cl0\cup\ll1)$}}%
    \put(20,10){\circle*{1.1}}%
    \put(21.5, 9){$\ll1=\rmin'(\ll0\cup\cl0\cup\ll1)$}%
    \put(20, 0){\circle*{1.1}}%
    \put(21, 0){$\emptyset=\rmin'(\Di_{\{0,1\}})$}%
    \put(20,20){\line(0, 1){10}}%
    \put(10,15){\line(2, 1){10}}%
    \put(10,15){\line(2,-1){10}}%
    \put(20, 0){\line(0, 1){10}}%
    \put(20,20){\line(2,-1){10}}%
    \put(20,10){\line(2, 1){10}}%
\end{picture}
\caption{The crystal lattice}\label{fig2}
\end{figure}
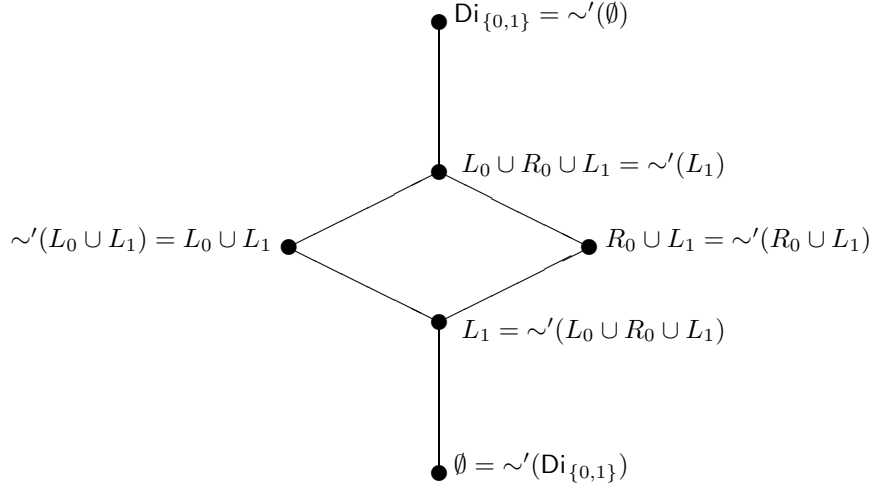
\begin{table}
\begin{gather*}  
  \begin{array}{|c|cccccc|}\hline
    \X\to'\Y&\Di_{\{0,1\}}&\ll0\cup\cl0\cup\ll1
    &\ll0\cup\ll1&\cl0\cup\ll1&\ll1&\emptyset\\\hline
    \Di_{\{0,1\}}&\Di_{\{0,1\}}
    &\emptyset&\emptyset&\emptyset&\emptyset&\emptyset\\
    \ll0\cup\cl0\cup\ll1&\Di_{\{0,1\}}&\ll1
    &\emptyset&\emptyset&\emptyset&\emptyset\\
    \ll0\cup\ll1&\Di_{\{0,1\}}&\ll0\cup\ll1
    &\ll0\cup\ll1&\emptyset&\emptyset&\emptyset\\
    \cl0\cup\ll1&\Di_{\{0,1\}}&\cl0\cup\ll1
    &\emptyset&\cl0\cup\ll1&\emptyset&\emptyset\\
    \ll1&\Di_{\{0,1\}}&\ll0\cup\cl0\cup\ll1
    &\ll0\cup\ll1&\cl0\cup\ll1&\ll1&\emptyset\\
    \emptyset&\Di_{\{0,1\}}&\Di_{\{0,1\}}
    &\Di_{\{0,1\}}&\Di_{\{0,1\}}
    &\Di_{\{0,1\}}&\Di_{\{0,1\}}\\\hline
  \end{array}
  \\
  \begin{array}{|c|cccc|}\hline
    \X|'\Y&\ll0&\cl0&\ll1&\cl1\\\hline
    \ll0&\ll0&\ll0\cup\cl0&\ll1&\cl1\\
    \cl0&\ll0\cup\cl0&\cl0&\ll1&\cl1\\
    \ll1&\ll1&\ll1&\ll1&\ll0\cup\cl0\cup\ll1\cup\cl1\\
    \cl1&\cl1&\cl1&\ll0\cup\cl0\cup\ll1\cup\cl1&\cl1\\\hline
  \end{array}
\end{gather*}  
\caption{Tables for the crystal lattice and $\mathfrak{S}_{\{0,1\}}$}
\label{crystaltables}
\end{table}
The crystal lattice first appears in Routley~\cite{MR544617}, where it
is attributed to R.\ K.\ Meyer; see \cite[pp.\ 65--6]{zbMATH03928963},
\cite[p.\ 250]{MR728950}, and \cite[pp.\ 95--7]{MR3728341}.  By
\cite[Theorems~9.8.1, 9.8.3]{MR3728341}, the crystal lattice
$\mathbf{Cr}$ is characteristic for the finitely axiomatized logic
$\CL$ \cite[p.\ 114]{MR3728341}.  We can obtain the crystal lattice
from $\sra_{\{0,1\}}$, which is isomorphic to the relation algebra
\alg{2}{83}, the second of 83 relation algebras listed in \cite[Ch.\
6, \SS\SS62--3]{MR2269199}.
\begin{theorem}
  The crystal lattice is isomorphic to a relativized subreduct of
  $\sra_{\{0,1\}}$.
\end{theorem}
\proof For a copy of the crystal lattice in $\sra_{\{0,1\}}$, let
\begin{align*}
  {Cr}& =\{\emptyset,\,\, \ll1,\,\, \ll0\cup\ll1,\,\,
  \cl0\cup\ll1,\,\, \ll0\cup\cl0\cup\ll1,\,\, \Di_{\{0,1\}}\},\\
  \mathbf{Cr}&=\<{Cr},\,\cup,\,\cap,\,\to',\,\rmin'\>.
\end{align*}
Inspection shows ${Cr}$ is closed under union, intersection,
relativized residuation, and relativized converse-complementation.
Comparison with \cite[p.\ 250]{MR728950} or \cite[pp.\
95--7]{MR3728341} shows $\mathbf{Cr}$ is the crystal lattice.  The
Hasse diagram and the action of $\rmin'$ are shown in Figure
\ref{fig2}, while $\to'$ is given in Table \ref{crystaltables}.
${Cr}$ is the union of two Sugihara chains of length 5 that intersect
in all but one relation. To get these two chains, delete either
$\ll0\cup\ll1$ or $\cl0\cup\ll1$ from ${Cr}$.  ${Cr}$ is also a set of
generators for $\sra_{\{0,1\}}$ (since conversion and complementation
are allowed).  Table \ref{crystaltables} shows the relativized
relative products of the diversity atoms of $\sra_{\{0,1\}}$.
\endproof
$\mathbf{Cr}$ is used in \cite[Theorem 3.22]{MR728950} for a proof of
the variable-sharing property that is simpler because it uses a
smaller algebra, with only six elements instead of eight, and the
2-element chains $\{<,\leq\}$ and $\{>,\geq\}$ in Belnap's proof are
replaced by singletons $\{\ll0\cup\ll1\}$ and $\{\cl0\cup\ll1\}$.
\subsection*{The Church lattice}
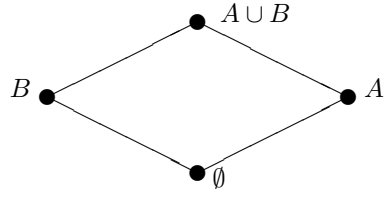
\begin{figure}
  \begin{picture}(35,14)(0,0)%
    \put(20  ,11 ){\circle*{1.1}}%
    \put(21.5,11 ){$\A\cup\B$}%
    \put(10  , 6 ){\circle*{1}}%
    \put( 7.5, 6 ){$\B$}%
    \put(30  , 6 ){\circle*{1}{$\A$}}%
    \put(20  , 1 ){\circle*{1.1}}%
    \put(21  , 0 ){$\emptyset$}%
    \put(10  , 6 ){\line(2, 1){10}}%
    \put(10  , 6 ){\line(2,-1){10}}%
    \put(20  ,11 ){\line(2,-1){10}}%
    \put(20  , 1 ){\line(2, 1){10}}%
\end{picture}
\caption{The Church diamond}\label{fig3}
\end{figure}
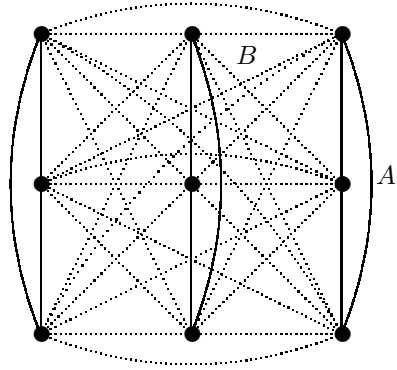
\begin{figure}
\setlength{\unitlength}{2mm}
\begin{picture}(22,22)(-1,-1)
\thinlines
\put( 0.1,20){\circle*{1}}
\put(10.1,20){\circle*{1}}
\put(20.1,20){\circle*{1}}
\put( 0.1,10){\circle*{1}}
\put(10.1,10){\circle*{1}}
\put(20.1,10){\circle*{1}}
\put( 0.1, 0){\circle*{1}}
\put(10.1, 0){\circle*{1}}
\put(20.1, 0){\circle*{1}}
\put(22.3,10){$A$}
\put(13,18){$B$}
\qbezier    ( 0, 0)( 0, 5)( 0,10)
\qbezier    ( 0, 0)(-4,10)( 0,20)
\qbezier[30]( 0, 0)( 5, 0)(10, 0)
\qbezier[42]( 0, 0)( 5, 5)(10,10)
\qbezier[60]( 0, 0)( 5,10)(10,20)
\qbezier[60]( 0, 0)(10,-4)(20, 0)
\qbezier[60]( 0, 0)(10, 5)(20,10)
\qbezier[60]( 0, 0)( 7,13)(20,20)
\qbezier    ( 0,10)( 0,15)( 0,20)
\qbezier[42]( 0,10)( 5, 5)(10, 0)
\qbezier[30]( 0,10)( 5,10)(10,10)
\qbezier[42]( 0,10)( 5,15)(10,20)
\qbezier[68]( 0,10)(10, 5)(20, 0)
\qbezier[60]( 0,10)(10,14)(20,10)
\qbezier[68]( 0,10)(10,15)(20,20)
\qbezier[60]( 0,20)( 5,10)(10, 0)
\qbezier[42]( 0,20)( 5,15)(10,10)
\qbezier[30]( 0,20)( 5,20)(10,20)
\qbezier[84]( 0,20)(13,13)(20, 0)
\qbezier[60]( 0,20)(10,15)(20,10)
\qbezier[60]( 0,20)(10,24)(20,20)
\qbezier    (10, 0)(10, 5)(10,10)
\qbezier    (10, 0)(14,10)(10,20)
\qbezier[30](10, 0)(15, 0)(20, 0)
\qbezier[42](10, 0)(15, 5)(20,10)
\qbezier[60](10, 0)(15,10)(20,20)
\qbezier    (10,10)(10,15)(10,20)
\qbezier[42](10,10)(15, 5)(20, 0)
\qbezier[30](10,10)(15,10)(20,10)
\qbezier[42](10,10)(15,15)(20,20)
\qbezier[60](10,20)(15,10)(20, 0)
\qbezier[42](10,20)(15,15)(20,10)
\qbezier[30](10,20)(15,20)(20,20)
\qbezier    (20, 0)(20, 5)(20,10)
\qbezier    (20, 0)(24,10)(20,20)
\qbezier    (20,10)(20,15)(20,20)
\end{picture}
\caption{Atoms of Church's relation algebra: $A=$ pairs connected by
  solid lines, $B=$ pairs connected by dotted lines}
\label{fig4}
\end{figure}
\begin{table}
  \begin{gather*}
    \begin{array}{|c|cccc|c|}\hline \X\to'\Y&
      \A\cup\B&\A&\B&\emptyset&\rmin'\\\hline
      \A\cup\B& \A\cup\B&\emptyset&\emptyset&\emptyset&\emptyset\\
      \A& \A\cup\B&\A&\B&\emptyset&\B\\
      \B& \A\cup\B&\emptyset&\A&\emptyset&\A\\
      \emptyset& \A\cup\B&\A\cup\B&\A\cup\B&\A\cup\B&\A\cup\B\\\hline
    \end{array}
    \\
    \begin{array}{|c|ccc|}\hline \X|\Y%
      &\Id&\A&\B\\\hline
      \Id&\Id&\A&\B\\
      \A&\A&\Id\cup\A&\B\\
      \B&\B&\B&\Id\cup\A\cup\B\\\hline
    \end{array}\\
  \end{gather*}
\caption{Tables for Church's relation algebra}
\label{churchtables}
\end{table}
The Church lattice \cite[p.\ 379]{MR728950} is also called Church's
diamond \cite[p.\ 277]{Schechter2005}.
\begin{theorem}\label{church}
  The Church lattice is the relativized reduct of a proper relation
  algebra on any set with 9 or more elements.
\end{theorem}
\proof On any set $\U$ with at least 9 elements, let $\V_1$, $\V_2$,
and $\V_3$ be a partition of $\U$ into pairwise disjoint sets, each
containing at least 3 elements. In the 9-element case, $\V_1$, $\V_2$,
and $\V_3$ are arranged in three columns as in Figure \ref{fig4}.  Let
\begin{align*}
  \U&=\V_1\cup\V_2\cup\V_3,\quad \Id=\{\<\u,\u\>\colon\u\in\U\}, \quad
  \Di=\{\<\u,\v\>\colon\u,\v\in\U,\,\u\neq\v\},\\
  \A&=\Di\cap\big((\V_1)^2\cup(\V_2)^2\cup(\V_3)^2\big),\quad
  \B=\bigcup\{\V_\i\times\V_\j:1\leq i,j\leq3,\,\i\neq\j\}.
\end{align*}
Then $\{\Id,\A,\B\}$ is a partition of $\U^2$ into relations that are
symmetric.  The eight unions of subsets of $\{\Id,\A,\B\}$ form a
proper relation algebra $\mathfrak{Ch}$ with $\{\Id,\A,\B\}$ as its
set of atoms. As noted in \SS2, $\mathfrak{Ch}$ is isomorphic to
relation algebra \alg47 \cite[Ch.\ 6, \SS56.13]{MR2269199}. The
relative products of atoms are shown in Table \ref{churchtables}.
Unions of symmetric relations are symmetric, so $\mathfrak{Ch}$ is a
symmetric proper relation algebra, called {\bf Church's relation
  algebra}.  Symmetric relation algebras are commutative because, by
axiom \thetag{r6} and the symmetry of both the factors and the
relative product, $\x\rp\y = \con{\x\rp\y} = \con\y\rp\con\x =
\y\rp\x$.  Not all symmetric relation algebras are dense, but
$\mathfrak{Ch}$ is dense.  The four diversity relations form Church's
diamond,
\begin{align*}
  {Ch}&=\{\A\cup\B,\A,\B,\emptyset\},\\
  \mathbf{Ch}&=\<{Ch},\,\cup,\,\cap,\,\to',\,\rmin'\>,
\end{align*}
with a Hasse diagram in Figure \ref{fig3}.  Tables for $\to'$ and
$\rmin'$ are in Table \ref{churchtables}.
\endproof
The Church lattice $\mathbf{Ch}$ validates the logic $\KR$, which is
axiomatized by axioms (R1)--(R13) in Table \ref{table2} along with
$(\X\land\rmin\X)\to\Y$.  The Lindenbaum algebra of $\KR$ is a
relation algebra \cite[Lemma 6.7]{MR2536403}.  (The method for
creating what is here called a ``Lindenbaum'' algebra is due to Tarski
\cite{zbMATH02525153}, \cite[Ch.  XII]{MR736686}; see
\cite{MR1117874}, \cite[p.\ 122, footnote 7]{MR6334}, \cite[p.\ 85,
footnote 4]{MR0124250}, \cite[p.\ 169, footnote 2]{MR781929}.)
$\mathbf{Ch}$ shows that $\KR$ is ``crypto-relevant'' \cite[p.\
379]{MR728950}, which means that the variable sharing property holds
for a formula $\X\to\Y$ if the only connective appearing is $\to$. To
show this, assign the variables in $\X$ to $\A\cup\B$ and the
variables in $\Y$ to $\A$.  Then $\X$ and $\Y$ are mapped to
$\A\cup\B$ and $\A$ since these are fixed by the operation $\to'$, but
$(\A\cup\B)\to'\A=\emptyset$ and the designated elements are
$\A\cup\B$ and $\A$, so $\X\to\Y$ is not valid in $\mathbf{Ch}$.
\subsection*{Meyer's RM84}
\begin{figure}
\begin{picture}(60,33)(-10,-1)%
  \put(20,30){\circle*{1}{$\{0,1,2,3,4,5,6\}$}}%
  \put(40,20){\circle*{1}{$\{0,3,5,6\}$}}%
  \put(20,20){\circle*{1}{$\{1,2,3,4,5,6\}$}}%
  \put(40,10){\circle*{1}{$\{3,5,6\}$}}%
  \put( 0,10){\circle*{1}}\put(-7,10){$\{1,2,4\}$}%
  \put(20,10){\circle*{1}}\put(21,9){$\{0\}$}%
  \put( 0,20){\circle*{1}}\put(-8.7,20){$\{0,1,2,4\}$}%
  \put(20, 0){\circle*{1}}\put(21.5,-1){$\emptyset$}%
  \put(20, 0){\line(2,1){20}}%
  \put( 0,10){\line(2,1){20}}%
  \put(20,10){\line(2,1){20}}%
  \put( 0,20){\line(2,1){20}}%
  \put( 0,10){\line(2,-1){20}}%
  \put(20,30){\line(2,-1){20}}%
  \put(20,20){\line(2,-1){20}}%
  \put( 0,20){\line(2,-1){20}}%
  \put( 0,10){\line(0, 1){10}}%
  \put(20,20){\line(0, 1){10}}%
  \put(40,10){\line(0, 1){10}}%
  \put( 20,0){\line(0, 1){10}}%
\end{picture}
\caption{Hasse diagram for RM84}
\label{fig5}
\end{figure}
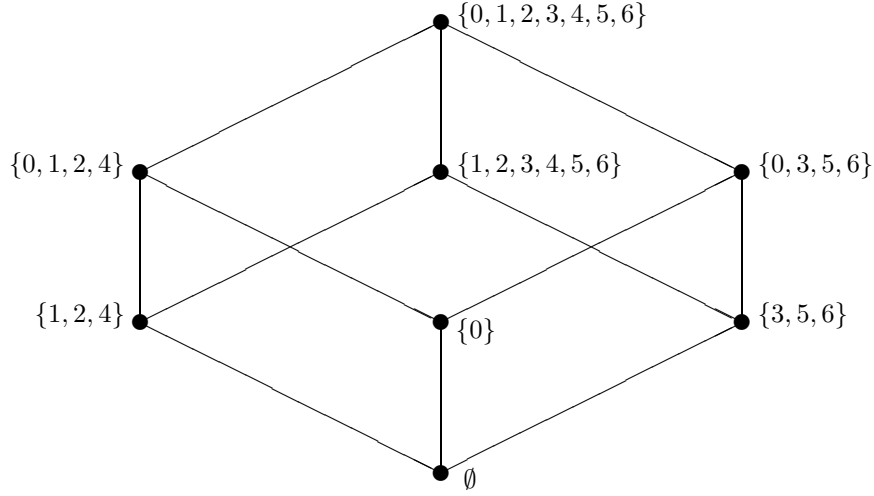
\begin{table}
\begin{gather*}
  \begin{array}{|c|c|}\hline%
    \X &\rmin\X \\\hline
    \emptyset &\U \\
    \{3,5,6\} &\{0,3,5,6\} \\
    \{1,2,4\} &\{0,1,2,4\} \\
    \D=\{1,2,3,4,5,6\} &\{0\} \\
    \{0\} &\D \\
    \{0,1,2,4\} &\{1,2,4\} \\
    \{0,3,5,6\} &\{3,5,6\} \\
    \U=\{0,1,2,3,4,5,6\} &\emptyset \\\hline
  \end{array}\\
  \begin{array}{|c|ccc|}\hline
        \X|\Y&\{0\}    &\{1,2,4\}&\{3,5,6\}\\\hline
        \{0\}&\{0\}    &\{1,2,4\}&\{3,5,6\}\\
    \{1,2,4\}&\{1,2,4\}&       \D&       \U\\
    \{3,5,6\}&\{3,5,6\}&       \U&       \D\\\hline
  \end{array}\\
  \begin{array}{|c|cccccccc|}\hline
    \X\to\Y&\3&\{3,5,6\}&\{1,2,4\}&\D&\{0\}&\{0,1,2,4\}&\{0,3,5,6\}&\U
    \\\hline\3&\U&\U&\U&\U&\U&\U&\U&\U
    \\\{3,5,6\}&\3&\{0\}&\3&\{0,3,5,6\}&\3&\3&\{0\}&\U
    \\\{1,2,4\}&\3&\3&\{0\}&\{0,1,2,4\}&\3&\{0\}&\3&\U
    \\\D&\3&\3&\3&\{0\}&\3&\3&\3&\U
    \\\{0\}&\3&\{3,5,6\}&\{1,2,4\}&\D&\{0\}&\{0,1,2,4\}&\{0,3,5,6\}&\U
    \\\{0,1,2,4\}&\3&\3&\3&\{1,2,4\}&\3&\{0\}&\3&\U
    \\\{0,3,5,6\}&\3&\3&\3&\{3,5,6\}&\3&\3&\{0\}&\U
    \\\U&\3&\3&\3&\3&\3&\3&\3&\U
    \\\hline
  \end{array}
\end{gather*}
\caption{Tables for RM84}
\label{rm84tables}
\end{table}
Anderson and Belnap \cite[p.\ 334]{MR0406756} present Meyer's lattice,
but they do not give it a name.  Instead, ``RM84'' is their name for
Meyer's theorem \cite[p.\,417]{MR0406756}, which says that if
$\X\to\Y$ is a theorem of $\RM$ then either $\X$ and $\Y$ share a
variable or both $\rmin\X$ and $\Y$ are theorems of $\RM$.  When
Routley, Plumwood, Meyer, and Brady \cite[p.\ 253]{MR728950} present
Meyer's lattice, they call it ``RM84'', as is done here.  In
\cite[Theorem 3.26]{MR728950} they show RM84 verifies all theorems of
$\RR$, but fails to satisfy any of eight particular formulas that
happen to be theorems of $\RM$.  The proper relation algebra
$\mathfrak{Rm}$, described here by subsets of the cyclic group of
order 7 and called {\bf Meyer's relation algebra}, is isomorphic to
relation algebra \alg33 \cite[Ch.\ 6, \SS58.8]{MR2269199}.
\begin{theorem}{\cite[Theorem 4.2]{MR2641636}}\label{meyer}
  {\normalshape RM84} is the relativized reduct the proper relation
  algebra $\mathfrak{Rm}$.
\end{theorem}
\proof Let $\U=\{0,1,2,3,4,5,6\}$, $\D=\{1,2,3,4,5,6\}$, and
\begin{align*}
  \Rm&=\{\U, \D, \{0,1,2,4\}, \{0,3,5,6\}, \{1,2,4\}, \{3,5,6\},
  \{0\}, \emptyset\}.
\end{align*}
We use $\Rm$ as an index set for eight binary relations on $\U$.  For
$\x,\y\in\U$, let $\x\equiv_7\y$ mean $\y-\x$ is divisible by $7$, and
for every $\X\subseteq\U$, define a relation on $\U$ by
\begin{align*}
  \rho(\X)&=\{\<\y,\z\>\colon\y,\z\in\U,\,\z+\x\equiv_7\y \text{ for
    some }\x\in\X\}.
\end{align*}
Then $\{\rho(\X)\colon\X\in\Rm\}$ is the universe of the proper
relation algebra $\mathfrak{Rm}$, which is an 8-relation subalgebra of
$\Re(\U)$.  Figure~\ref{fig5} shows the Hasse diagram for sets in
$\Rm$ and their images under $\rho$.  The images of $\{0\}$,
$\{1,2,4\}$, and $\{3,5,6\}$ are atoms of $\mathfrak{Rm}$.  The
converse-complements and relative products in Table \ref{rm84tables}
are stated in terms of sets in $\Rm$. The entry for $\X,\Y\in\Rm$ is
the set $\Z\in\Rm$ such that $\rho(\Z)=\rho(\X)|\rho(\Y)$.
Converse-complements and relative products can also be computed
directly by the rules $\rmin\X=\{0-_7\x:\x\notin\X,\,\x\in\U\}$ and
$\X|\Y = \{\x+_7\y \colon \x\in\X,\y\in\Y\}$, where $-_7$ and $+_7$
are subtraction and addition {\it modulo} 7.  As described here, RM84
is the direct reduct of $\mathfrak{Rm}$.
\endproof
\section{A relational completeness %
  theorem for $\RR$-mingle}\label{sect9}
\begin{table}
\begin{align*}
  &\tag{R1}\A\to\A\\
  &\tag{R2}(\A\to\B)\to((\B\to\C)\to(\A\to\C))\\
  &\tag{R3}\A\to((\A\to\B)\to\B)\\
  &\tag{R4}(\A\to(\A\to\B))\to(\A\to\B)\\
  &\tag{R5}(\A\land\B)\to\A\\
  &\tag{R6}(\A\land\B)\to\B\\
  &\tag{R7}((\A\to\B)\land(\A\to\C))\to(\A\to(\B\land\C))\\
  &\tag{R8}\A\to(\A\lor\B)\\
  &\tag{R9}\B\to(\A\lor\B)\\
  &\tag{R10}((\A\to\C)\land(\B\to\C))\to((\A\lor\B)\to\C)\\
  &\tag{R11}(\A\land(\B\lor\C))\to((\A\land\B)\lor\C)\\
  &\tag{R12}(\A\to\rmin\B)\to(\B\to\rmin\A)\\
  &\tag{R13}\rmin\rmin\A\to\A\\
  &\tag{R14}\A\to(\A\to\A)
\end{align*}
\caption{Axioms of $\RM$}\label{table2}
\end{table}
The logic $\RR$-mingle, or $\RM$, was created by Dunn and McCall from
Anderson and Belnap's relevance logic $\RR$ by adding the mingle axiom
$\A\to(\A\to\A)$; see \cite[\SS8.15, \SS27.1.1] {MR0406756}.  The
rules of deduction for both $\RR$ and $\RM$ are {\it Adjunction}
(infer $\A\land\B$ from $\A$ and $\B$) and {\it modus ponens} (infer
$\B$ from $\A\to\B$ and $\A$).  An axiom set for $\RM$ is shown in
Table \ref{table2}; see \cite[p.\ 341]{MR0406756} or \cite[pp.\
xxiii--xxvi]{MR1223997}.

If $\Su$ is a Sugihara chain and the connectives of $\RM$ are
interpreted as the corresponding operations (with the same names) in
$\Su$, then any function from the propositional variables to elements
of $\Su$ extends uniquely to a homomorphism from the algebra of
formulas to elements of $\Su$.  A formula is {\bf valid in} $\Su$ if
it is sent to a designated element by every such homomorphism.  Meyer
\cite[pp.\ 413--4, Corollaries 3.1, 3.5]{MR0406756} proved that the
theorems of $\RM$ are the formulas valid in all finite Sugihara
chains, and that the theorems of $\RM$ are the formulas valid in
$\Su_{\Zn^*}$.  These results, together with Theorem~\ref{iso}, imply
that $\RM$ is complete with respect to the following class of
algebras.
\begin{definition}\label{KRM}
  Let $\Ka=\<\K,\cup,\cap,\to,\rmin\>$, where
  \begin{enumerate}
  \item[(k1)] $\K$ is a non-empty set of binary relations on a set
    $\U$, called the {\bf base} of $\Ka$,
  \item[(k2)] $\K$ is closed under the operations $\cup$, $\cap$,
    $\to$, and $\rmin$, defined in Table \ref{table1} using the base
    $\U$.
  \end{enumerate}
  A formula $\A$ is {\bf valid in} the algebra $\Ka$ if, for every
  homomorphism $\h$ from the algebra of formulas to $\Ka$, $\h(\A)$
  contains the identity relation on the base of $\Ka$.  Let $\KRM$ be
  the class of algebras $\Ka=\<\K,\cup,\cap,\to,\rmin\>$ such that
  \thetag{k1}, \thetag{k2}, and
  \begin{enumerate}
  \item[(k3)] $\A|\B=\B|\A$ for all $\A,\B\in\K$,
  \item[(k4)] $\A\subseteq\A|\A$ for all $\A\in\K$,
  \item[(k5)] $\A|\A\subseteq\A$ for all $\A\in\K$.
  \end{enumerate}
  A formula is {\bf valid in }$\KRM$ if it is valid in every algebra
  in $\KRM$.
\end{definition}
Since $\K$ is not empty by \thetag{k1}, the algebra $\Ka$ determines
the base according to the formula $\U=\{\a:\<\a,\a\>\in\A
\cup\rmin\A,\,\A\in\K\}$. Condition \thetag{k2} implies that $\Ka$ is
also closed under $|$, since $\A|\B=\rmin(\B\to\rmin\A)$ by
\eqref{fuse}. From Theorem~\ref{iso} and Meyer's results we get the
following completeness theorem.
\begin{theorem}[{\cite[Theorem~6.2(iii)]{MR2641636}}]\label{comp}
  The theorems of $\RM$ are the formulas valid in $\KRM$.
\end{theorem}\proof
The axioms of $\RM$ are valid in $\KRM$ by Theorem \ref{meanings}
below.  Validity is preserved by Adjunction, for if $\Id\subseteq\A$
and $\Id\subseteq\B$ then $\Id\subseteq\A\cap\B$, and validity is
preserved by {\it modus ponens}, for if $\Id\subseteq\A\to\B$ and
$\Id\subseteq\A$, then $\Id=\Id|\Id\subseteq\A|(\A\to\B)\subseteq\B$
by Lemma \ref{recall}\ref{little} below.  Therefore all theorems of
$\RM$ are valid in $\KRM$.  For the converse, suppose $\X$ is not a
theorem of $\RM$.  By Meyer's theorems, $\X$ fails in $\Su_{\Zn^*}$,
$\Su_{\Zn+\Zn}$, and in every sufficiently large finite normal
Sugihara chain.  By Theorem~\ref{iso}, these algebras are isomorphic
to algebras in $\KRM$. Therefore there are algebras in $\KRM$ in which
$\X$ fails to be valid.
\endproof
We assume for the rest of this section that $\U$ is a set and
$\Ka=\<\K,\cup,\cap,\to,\rmin\>$ is an algebra satisfying conditions
(k1) and (k2) of Definition \ref{KRM}.  In each formula, we interpret
the connectives $\lor$, $\land$, $\to$, and $\rmin$ as the operations
$\cup$, $\cap$, $\to$, and $\rmin$, respectively. Thus every formula
denotes a relation that depends on the interpretation of its
variables. A formula is valid in $\Ka$ if it denotes a relation that
contains $\Id$, the identity relation on $\U$, no matter how its
variables are interpreted.  Implications are analyzed as inclusions
because of the following lemma.
\begin{lemma}[{\cite[Theorem~5.1(17)]{MR2641636}}]\label{inc}
  For all $\A,\B\subseteq\U^2$, $\Id\subseteq\A\to\B$ iff
  $\A\subseteq\B$.
\end{lemma}
According to Lemma~\ref{inc}, the validity of each axiom of $\RM$ can
be equivalently expressed as an inclusion between binary relations.
For example, (R1) is valid in $\Ka$ just because the inclusion
$\A\subseteq\A$ always holds.  Evidently (R1), (R5), (R6), (R7), (R8),
(R9), (R10), (R11), and (R13) are true under the set-theoretical
meanings assigned to the connectives, by \cite[Theorem~5.1(32), (33),
(34), (35), (36), (37), (38), (39), (40)]{MR2641636}, respectively. To
analyze the remaining axioms of $\RM$ we recall some other results
from \cite{MR2641636}.
\begin{lemma}{\cite[Theorem~5.1(18)(19)(21)(22)]{MR2641636}}
  \label{recall} For all $\A,\B,\C\subseteq\U^2$,
  \begin{enumerate}
  \item \label{above} $\A\to(\B\to\C)=(\B|\A)\to\C$,
  \item \label{little} $\A|(\A\to\B)\subseteq\B$,
  \item \label{mono} if $\A\subseteq\B$ then
    $\B\to\C\subseteq\A\to\C$ and $\C\to\A\subseteq\C\to\B$.
  \end{enumerate}
\end{lemma}
Axioms (R1), (R5)--(R11), and (R13) are valid for all binary
relations.  The remaining five axioms do not hold for all relations,
but will hold under conditions on the relations that occur in them,
and in some cases are equivalent to those conditions.  We now analyze
(R2), (R3), (R4), (R12), and (R14).  By \cite[Theorem
5.1(55)]{MR2641636}, (R2) holds whenever $\B\to\C$ and $\A\to\B$
commute, but (R2) also holds under the weaker hypothesis of
Lemma~\ref{R2} below, because inclusion in only one direction is
needed.  (R2) holds if $\Ka$ is commutative under relative
multiplication, but fails in some non-commutative examples that have
16 relations. On the other hand, (R2) is valid when recast as a rule
of inference, for if $\A\to\B$ contains the identity relation then so
does $(\B\to\C)\to(\A\to\C)$ \cite[Theorem 5.1(29)]{MR2641636}. This
also follows immediately from Lemma~\ref{inc} and Lemma
\ref{recall}\ref{mono}.
\begin{lemma}\label{R2}
  For all $\A,\B,\C\subseteq\U^2$, if $(\B\to\C)|(\A\to\B) \subseteq
  (\A\to\B)|(\B\to\C)$ then $\A\to\B \subseteq (\B\to\C)\to(\A\to\C)$
  and \thetag{R2} is valid, but the converse may fail.
\end{lemma}\proof First we prove the assumption implies the validity
of \thetag{R2}.
\begin{align*}
  &(\A|(\B\to\C))|(\A\to\B)\\
  &=\A|((\B\to\C)|(\A\to\B))
  &&\text{$|$ is associative}\\
  &\subseteq\A|((\A\to\B)|(\B\to\C))
  &&\text{assumption, $|$ is monotonic}\\
  &=(\A|(\A\to\B))|(\B\to\C)
  &&\text{$|$ is associative}\\
  &\subseteq\B|(\B\to\C)
  &&\text{Lemma \ref{recall}\ref{little}, $|$ is monotonic}\\
  &\subseteq\C,&&\text{Lemma \ref{recall}\ref{little}}
\end{align*}
hence
\begin{align*}
  \C\to\C&\subseteq((\A|(\B\to\C))|(\A\to\B))\to\C
  &&\text{Lemma \ref{recall}\ref{mono}}\\
  &=(\A\to\B)\to((\A|(\B\to\C))\to\C)
  &&\text{Lemma \ref{recall}\ref{above}}\\
  &=(\A\to\B)\to((\B\to\C)\to(\A\to\C)). %
  &&\text{Lemma \ref{recall}\ref{above}}
\end{align*}
Since $\Id\subseteq\C\to\C$, it follows that \thetag{R2} is valid,
\ie,
\begin{align*}
  \Id&\subseteq(\A\to\B)\to((\B\to\C)\to(\A\to\C)).
\end{align*}
By Lemma~\ref{inc}, $\A\to\B\subseteq(\B\to\C)\to(\A\to\C)$. However,
the assumption and this conclusion are not equivalent. To see this,
let $\C=\U^2$. Then $\A\to\C=\U^2$, hence
\begin{align*}
  (\B\to\C)\to(\A\to\C)&=(\B\to\C)\to\U^2=\U^2,
\end{align*}
so the conclusion holds for all $\A$ and $\B$. Since $\A\to\C=\U^2$,
the assumption becomes $\U^2|\left(\A\to\B\right) \subseteq
\left(\A\to\B\right)|\U^2$.  But this inclusion will fail whenever
$\A\to\B$ is not empty and has a domain that is not all of $\U$.
\endproof
By \cite[Theorem~5.1(54)]{MR2641636}, (R3) holds whenever $\A$
and $\A\to\B$ commute. In fact, it holds under a weaker hypothesis to
which it is not equivalent.
\begin{lemma}\label{R3}
  If $(\A\to\B)|\A\subseteq\A|(\A\to\B)$ then \thetag{R3} is valid.
  The converse may fail.
\end{lemma}
\proof We have $(\A\to\B)|\A\subseteq\B$ by the hypothesis and Lemma
\ref{recall}\ref{little}. Let $\C=\A\to\B$, so that
$\C|\A\subseteq\B$. This formula can be rewritten as
$\min\B\cap(\C|\A)=\emptyset$.  This is equivalent to
$(\conv\C|\min\B)\cap\A=\emptyset$, which is in turn equivalent to
$\A\subseteq\min{\conv\C|\min\B}$, but $\min{\conv\C|\min\B}=\C\to\B$,
so $\A\subseteq\C\to\B$. Hence $\A\subseteq(\A\to\B)\to\B$ and
\thetag{R3} is valid by Lemma \ref{inc}.  This conclusion does not
imply the hypothesis, for if $\B=\U^2$, then $\A\to\B=
\min{\conv\A|\min\B}= \min{\conv\A|\min{\U^2}}= \U^2$, so the
hypothesis is equivalent to $\U^2|\A\subseteq\A|\U^2$, which fails if
$\A$ is a relation on $\U$ whose domain is not all of $\U$. On the
other hand, the conclusion of Lemma~\ref{R3} holds since
$(\A\to\B)\to\B=\min{\conv{\U^2}|\min{\U^2}}=\U^2$.
\endproof
\begin{lemma}\label{R4}
  $\thetag{R4}$ is valid if and only if $\A$ is dense.
\end{lemma}
\proof By \cite[Theorem~5.1(56)]{MR2641636}, (R4) is valid whenever
$\A$ is a dense relation, for if $\A\subseteq\A|\A$ then
$\A\to(\A\to\B)\subseteq\A\to\B$.  Suppose (R4) is valid when
$\B=\rmin\Id$.  Since $\A\to\rmin\Id = \min{\conv\A|\min{\rmin\Id}} =
\min{\conv\A|\conv\Id} = \rmin\A$, (R4) is equivalent to
$(\A\to\rmin\A) \to \rmin\A$, which is valid if and only if
$\A\to\rmin\A\subseteq\rmin\A$, by Lemma \ref{inc}. This last
inclusion can be equivalently transformed by the definitions of
$\rmin$ and $\to$ first into $\min{\conv\A|\min{\rmin\A}} \subseteq
\min{\conv\A}$, then $\conv\A \subseteq \conv\A|\conv\A$, and finally
$\A\subseteq\A|\A$, which asserts that $\A$ is dense.
\endproof
The contraposition axiom (R12) is valid whenever $\A|\B=\B|\A$ by
\cite[Theorem~5.1(53)]{MR2641636}, but it is actually equivalent to
$\A|\B\subseteq\B|\A$.
\begin{lemma}\label{R12}
\thetag{R12} is valid if and only if $\A|\B\subseteq\B|\A$.
\end{lemma}
\proof By Lemma \ref{inc}, \thetag{R12} is valid if and only if
$\A\to\rmin\B\subseteq \B\to\rmin\A$ for all $\A,\B\subseteq\U^2$.
Since $\A\to\rmin\B = \min{\conv\A|\min{\rmin\B}} =
\min{\conv\A|\conv\B}$ and $\B\to\rmin\A = \min{\conv\B|\conv\A}$,
this inclusion is equivalent to $\conv\B|\conv\A \subseteq
\conv\A|\conv\B$.  Taking converses of both sides, we get the
equivalent inclusion $\A|\B \subseteq \B|\A$.
\endproof
By \cite[Theorem~5.1(63)]{MR2641636}, (R14) holds if $\A$ is a
transitive relation, but (R14) is actually equivalent to the
transitivity of $\A$.
\begin{lemma}\label{R14}
  \thetag{R14} is valid if and only if $\A$ is transitive.
\end{lemma}\proof
By Lemma \ref{inc}, \thetag{R14} is valid if and only if
$\A\subseteq\A\to\A$.  This inclusion can be equivalently restated
first as $\A \subseteq \min{\conv\A|\min\A}$, then
$\A\cap(\conv\A|\min\A) = \emptyset$, then $\A|\A\cap\min\A =
\emptyset$, and finally $\A|\A\subseteq\A$, which asserts that $\A$ is
transitive.
\endproof
The following theorem gathers together the observations above and
confirms that the axioms of $\RM$ are valid in $\KRM$, completing the
proof of Theorem \ref{comp}.
\begin{theorem}\label{meanings}
  Let $\Ka=\<\K,\cup,\cap,\to,\rmin\>$ be an algebra satisfying
  conditions \thetag{k1} and \thetag{k2} in Definition \ref{KRM}. Then
  \begin{enumerate}
  \item\label{1} \thetag{R1}, \thetag{R5}, \thetag{R6}, \thetag{R7},
    \thetag{R8}, \thetag{R9}, \thetag{R10}, \thetag{R11}, \thetag{R13}
    are valid in $\Ka$,
  \item\label{2} \thetag{R2} and \thetag{R3} are valid in $\Ka$ if
    \thetag{k3}, but neither is equivalent to \thetag{k3},
  \item\label{3} \thetag{R4} is valid in $\Ka$ if and only if
    \thetag{k4},
  \item\label{4} \thetag{R12} is valid in $\Ka$ if and only if
    \thetag{k3},
  \item\label{5} \thetag{R14} is valid in $\Ka$ if and only if
    \thetag{k5}.
  \end{enumerate}
\end{theorem}
\section{Interpreting formulas as relations}
\label{sect10}
Theorems \ref{comp} and \ref{meanings} suggest an alternative approach
to $\RM$.  Instead of adopting 14 axioms and two rules, simply define
$\RM$ as the set of formulas valid in $\KRM$.  It is then a theorem
that $\RM$ can be axiomatized by \thetag{R1}--\thetag{R14} and the
rules of Adjunction and {\it modus ponens}.

To explore the theorems and rules of $\RM$, assume that $\Ka = \<\K,
\cup, \cap, \to,\rmin\>$ is an algebra satisfying \thetag{k1} and
\thetag{k2} in Definition \ref{KRM}.  Even if $\Ka$ does not satisfy
\thetag{k3}, \thetag{k4}, or \thetag{k5}, nine of the axioms of $\RM$
are valid in $\Ka$ by Theorem \ref{meanings}\ref{1}.  Since
$\A\cup\rmin\A$ always contains the identity relation on $\U$,
$\A\lor\rmin\A$ is also valid in $\Ka$.  Thus $\A\lor\rmin\A$ is a
theorem of $\RM$.

Simple counterexamples show that $\B\to(\A\lor\rmin\A)$ and
$(\A\land\rmin\A)\to\C$ need not be valid in $\Ka$ and are not
theorems of $\RM$.  Counterexamples for relations in general can be
found on a 2-element set, but for $\RM$ it is more appropriate to use
the Sugihara chain $\{\emptyset,\,<,\,\leq,\,\Q^2\}$ from Table
\ref{ten}, which is an algebra in $\KRM$.  Just let $\A$, $\B$, and
$\C$ be $<$, $\Q^2$, and $\emptyset$, respectively.  The formula
$((\A\lor\C) \land \rmin\A) \to\C$ expressing Extensional Disjunctive
Syllogism also fails under the same assignment, so it is also not a
theorem of $\RM$.

The proof of Theorem \ref{comp} shows the rules of Adjunction and {\it
  modus ponens} preserve validity. Unlike the corresponding axiom,
Extensional Disjunctive Syllogism (to infer $\B$ from $\A\lor\B$ and
$\rmin\A$) is admissible.  To show this, assume $\Id\subseteq\A\cup\B$
and $\Id\subseteq\rmin\A$.  The second hypothesis is equivalent to
$\Id \subseteq \min{\conv\A}$. Taking the converse of both sides, we
get $\Id \subseteq \min\A$. By the first hypothesis, $\Id \subseteq
\min\A \cap (\A\cup\B) = \min\A\cap\B \subseteq \B$.  Intensional
Disjunctive Syllogism is to infer $\B$ from $\A+\B$ and $\rmin\A$,
where $\A+\B$ is intensional disjunction.  Since $\A+\B$ is defined as
$\rmin\A\to\B$ \cite[\SS27.1.4]{MR0406756}, this rule is an instance
of {\it modus ponens}.  The E-rule~\cite[p.\ 8]{MR3728341}, also
called BR1~\cite[p.\ 289]{MR728950} and R5~\cite[p.\ 193]{MR3728341},
is to infer $(\A\to\B)\to\B$ from $\A$.  Assume $\Id\subseteq\A$. Then
$\min\B = \conv\Id|\min\B \subseteq \conv\A|\min\B = \min{\A\to\B}$,
hence $\A\to\B\subseteq\B$, so $\Id\subseteq(\A\to\B)\to\B$ by Lemma
\ref{inc}. The admissibility of Suffixing, Contraposition, and several
other rules can be proved similarly at this point, without any appeal
to commutativity, density, or transitivity.

To achieve $\RM$, assume that $\Ka$ also satisfies \thetag{k3},
\thetag{k4}, and \thetag{k5}, so that $\Ka\in\KRM$. Then axioms
\thetag{R2}, \thetag{R3}, \thetag{R4}, \thetag{R12}, and \thetag{R14}
are also valid in $\Ka$ by Theorem
\ref{meanings}\ref{2}\ref{3}\ref{4}\ref{5}.  Alternate proofs of the
admissibility of various rules, such as the E-rule, Contraposition,
and Suffixing, are possible using commutativity.  A significant
example of a theorem of $\RM$ that requires all three hypotheses of
commutativity, density, and transitivity is $(\A\to\B)\lor(\B\to\A)$.
Meyer called this formula RM64, ``Simple order''. Discussing its
significance, he wrote \cite[pp.\ 397--8]{MR0406756},
\begin{quote}
  ``RM63 and RM64, in fact, decide that $\RM$ represents a much longer
  step in the direction of classical logic (and, for that matter, in
  the direction of an extensional approach to sentential logic) than
  one would have thought from the heuristic considerations by which we
  motivated its axioms and rules. \dots\par
  ``RM64 leaves shattered in the dust much of the motivation to which
  previous opponents of the paradoxes have appealed. But this just
  goes to show that one can have many reasons for disliking the
  paradoxes; one very plausible ground for disliking them is that they
  turn every minor inconsistency into a catastrophe. From this charge,
  $\RM$ is yet free. If in other respects it moves in the direction of
  classical logic, there is as yet no reason to rue that fact.''
\end{quote}
The following lemma gives a relational proof that
$(\A\to\B)\lor(\B\to\A)$ is a theorem of $\RM$.
\begin{lemma}
  If $\<\K,\cup,\cap,\to,\rmin\>\in\KRM$ then $\Id \subseteq (\A\to\B)
  \cup (\B\to\A)$ for all $\A,\B\in\K$.
\end{lemma}
\proof Note that $\K$ is closed under $|$ by \eqref{fuse}.  Assume
$\A,\B\in\K$, and let $\C = (\A\to\A) \cap(\B\to\B).$ Apply $\rmin$ to
both sides and use \eqref{fuse} to get 
\begin{equation*}
  \rmin\C = \rmin(\A\to\A) \cup
  \rmin(\B\to\B) = (\rmin\A|\A) \cup (\rmin\B|\B).
\end{equation*}
It follows that $\rmin\A|\A \subseteq \rmin\C$ and $\rmin\B|\B
\subseteq\rmin\C$, so $(\rmin\A|\A) | (\rmin\B|\B)
\subseteq\rmin\C|\rmin\C$ by the monotonicity of $|$.  By the
associativity of $|$ and our assumption that $|$ is commutative on
relations in $\K$, $(\rmin\A|\B) | (\rmin\B|\A) \subseteq
\rmin\C|\rmin\C.$ Let $\D=(\rmin\A|\B) \cap (\rmin\B|\A).$ From
$\D,\rmin\C\in\K$ it follows that $\D$ is dense and $\rmin\C$ is
transitive, so by the monotonicity of $|$,
\begin{equation*}
  \D\subseteq\D|\D\subseteq (\rmin\A|\B)| (\rmin\B|\A)
  \subseteq\rmin\C|\rmin\C\subseteq\rmin\C.
\end{equation*}
By applying $\rmin$ to both sides, \eqref{fuse}, and Lemma \ref{inc},
we conclude that
\begin{equation*}
  \Id\subseteq\C\subseteq\rmin\D=(\B\to\A)\cup(\A\to\B).
\end{equation*}
\endproof
RM64 is one reason given by Anderson and Belnap for the title of
\cite[\SS29.5]{MR0406756}, ``Why we don't like mingle.'' They describe
how to prove RM64, using axioms and rules, from the ``unhappy
theorem'' $\A\to(\rmin\A\to\A)$. Their suggestions for proving the
latter formula include the mingle axiom, contraposition, and
permutation, thereby invoking both transitivity and commutativity, but
commutativity is not needed. By Lemma \ref{inc}, $\A\to(\rmin\A\to\A)$
is valid if $\A\subseteq\rmin\A\to\A$ for every relation $\A\in\K$.
Since $\rmin\A\to\A=\rmin(\rmin\A|\rmin\A)$, this inclusion is
equivalent to $\rmin\A|\rmin\A\subseteq\rmin\A$, which asserts that
$\rmin\A$ is transitive, as is indeed the case for every $\A\in\K$.

The Routley-Meyer semantics are called relational because every
relevant model structure contains a ternary relation.  Instead of a
ternary relation, the $\RM$ model structures of Dunn
\cite{zbMATH03513746} have a binary accessibility relation, which
corresponds to the inclusion relation in the Sugihara chain $\Ch_\I$
\cite[\SS7]{zbMATH03513746}.  In both cases, the elements of the model
structures are objects with no further structure.  Formulas are
interpreted as sets of unstructured objects.  This feature is
advantageous because it provides more general interpretations than
semantics with special objects.  For example, Dunn \cite{MR3297467}
asks,
\begin{quote}
  ``What could be more natural than to interpret $Rabc$ as that in the
  context of the information $a$, the information $b$ is relevant to
  the information $c$?''  
\end{quote}
Three more examples (with variations) are presented by the eleven
authors of \cite{MR2914447}, based on three ways of grouping the
arguments of $R$, called Modal (Absence-of-Counterexample)
Conditionals: $Rx\<yz\>$, Conditionals as Operators: $R\<xy\>z$, and
Conditional Logics: $Rx\>y\<z$. These interpretations address some
issues, explained on \cite[p.\ 599]{MR2914447}.
\begin{quote}
  ``The story goes like this: whereas the binary relation invoked by
  Kripke in the semantics of modal logics has several philosophically
  interesting and revealing interpretations (as relative possibility,
  or as a temporal ordering, or as the relation of
  being-morally-ideal-from-the-point-of-view-of, or \dots), the
  ternary relation invoked by Routley and Meyer has no such standardly
  accepted interpretations/applications. `Sure,' the objector says (it
  helps here to imagine the hint of a sneer), `there are mathematical
  structures of the sort described by Routley and Meyer, and those
  structures bear important and interesting relations to the logics
  described by Anderson and Belnap, but these logics were supposed to
  tell us something interesting about \emph{conditionality}, or at
  least some important kind of conditionality, and it would take more
  than just abstract mathematical structures to tell us \emph{that}. I
  want to know what it is that \emph{instantiates} these structures
  that has anything to do with conditionals.' ''
\end{quote}
Although ``we say nothing about negation'' and ``this paper isn't
about negation'' \cite[footnote 4]{MR2914447}, elsewhere Dunn
\cite{MR3297467} observed,
\begin{quote}
  ``The `Routley-Star' has come under a lot of criticism both from
  those within and outside of the relevance logic community, and was
  more of a focus of Copeland's \cite{MR551278} critical review than
  the ternary accessibility relation.''
\end{quote}
This is reflected by van Benthem \cite{vanBenthem1984}, in his review
of Copeland \cite{MR551278}.
\begin{quote}
  ``Relevance logic is a subject whose motivation has turned out to be
  surprisingly difficult to capture in an enlightening and convincing
  semantics.  \dots the only general approach to date is a rather
  abstract possible-worlds framework, proposed by R.\ and V.\ Routley.
  Here, relevant implication is explicated through some ternary
  `perspective' relation among worlds, while the account of negation
  employs an additional `reversal' operation upon worlds.  It fails to
  satisfy the `requirements which distinguish an illuminating and
  philosophically significant semantics from a merely formal model
  theory.' ''
\end{quote}
Much later in the review, van Benthem observes,
\begin{quote}
  ``Postulated operations in model structures may be viewed as
  theoretical terms, truth-definitions rather as some kind of
  correspondence principles.  Any demand for `realism' ought to take
  these different roles into account.\par
  Even in this more charitable perspective, the Routley semantics
  still has to prove its mettle. On the realistic side, its model
  structures ought to admit of, if not a natural linguistic anchoring,
  then at least one mathematical `standard example', providing some
  food for independent reflection.''
\end{quote}
In our analysis, formulas are interpreted as sets of objects that
\emph{do} have structure. Each object is a binary relation.  Theorem
\ref{comp} says that a formula is a theorem of $\RM$ if and only if it
contains the identity relation when regarded as a relation belonging
to a set of transitive dense binary relations that commute under
relative multiplication. We call this the formulas as relations
approach.

Suppose $\gc\A$ is an atomic proper relation algebra on a set $\U$.
Let $\gc\U=\<\at,\R,\conv{},\I\>$ be the atom structure of $\gc\A$.
If $\gc\A$ happens to be commutative and dense, then $\gc\U$ is a
relevant model structure. The commutative dense atomic proper relation
algebras discussed in this paper are Belnap's relation algebra
$\mathfrak\M_0$ in \SS\ref{sect2}, Sugihara's relation algebra
$\sra_\Zn$ in \SS\ref{sect4} and, more generally $\sra_\I$ for any
$\I\subseteq\Zn$ in Theorem \ref{lemma1}, Church's relation algebra
$\mathfrak{Ch}$ in Theorem \ref{church}, and Meyer's relation algebra
$\mathfrak{Rm}$ in Theorem \ref{meyer}.  In any case, the atoms in
$\at$ are relations on $\U$, the unary operation $\conv{}$ of $\gc\U$
produces the converse of each atom, the distinguished element $\I$ of
$\gc\U$ is the set of atoms contained in the identity relation of
$\U$, and the ternary relation $\R$ of $\gc\U$ is the set of triples
$\<\a,\b,\c\>$ of relations that are atoms of $\gc\A$ and satisfy
$\a|\b\supseteq\c$.  In brief, the ternary relation is
set-theoretically defined as $\a|\b\supseteq\c$ and the Routley star
${}^*$ is conversion $\conv{}$. These interpretations of $\R$ and
${}^*$ are complete for $\RM$ according to Theorem \ref{comp}. They
may satisfy the objector described in \cite{MR2914447}, who wants ``to
know what it is that \emph{instantiates} these structures''.  For
$\RM$ at least, we could echo Dunn and ask, what could be more natural
than to interpret $Rabc$ as $\a|\b\supseteq\c$?  For van Benthem and
Copeland we suggest interpreting the ``additional `reversal'
operation'' as conversion.  In the light of Theorem \ref{comp}, $\RM$
could serve as a mathematical ``standard example'' sought by van
Benthem.
\section{Summary and problems}\label{sect11}
Table \ref{defred} summarizes our results that some finite lattices
with operators and all subalgebras of three countably infinite
Sugihara chains are isomorphic to definitional subreducts of proper
relation algebras.  The widespread occurrence of representability
where it was not previously suspected could lead to further thoughts
of a fundamental, perhaps philosophical nature. This idea is
elaborated in Problem \ref{prob4}, and there is further mathematical
work proposed in Problems \ref{prob1}, \ref{prob2}, and \ref{prob3}.
\begin{table}
  \begin{tabular}{|lllll|}\hline
    Name&Method&Type&PRA &Full PRA\\\hline
    Belnap $\MM$
    &Direct&reduct&$\sra_{\{0\}}$&$\Re(\U_{\{0\}})$\\
    Sugihara $\Su_{\Zn+\Zn}$, $\Su_{\Zn^*}$, finite even
    &Direct&subreduct&$\sra_\Zn$&$\Re(\U_\Zn)$\\
    Sugihara $\Su_{\Zn+\Zn}$, $\Su_\Zn$, $\Su_{\Zn^*}$, finite
    &Relativized&subreduct&$\sra_\Zn$&$\Re(\U_\Zn)$\\
    crystal $\mathbf{Cr}$
    &Relativized&subreduct&$\sra_{\{0,1\}}$&$\Re(\U_{\{0,1\}})$\\
    Church $\mathbf{Ch}$
    &Relativized&reduct&$\mathfrak{Ch}$&$\Re(9)$\\
    Meyer RM84&Direct&reduct&$\mathfrak{Rm}$&$\Re(7)$\\\hline
  \end{tabular}
  \caption{Definitional subreducts of proper relation algebras (PRAs)}
  \label{defred}
\end{table}
\begin{problem}\label{prob1}
  What other algebras in the relevance logic literature are isomorphic
  to definitional subreducts of proper relation algebras? For example,
  is every uncountable Sugihara chain isomorphic to a definitional
  subreduct of a proper relation algebra?
\end{problem}
\subsection*{Remarks on Problem \ref{prob1}}
Each of the logics $\CL$, $\BM$, and $\RM$ is characterized by a
single lattice.
\begin{itemize}
\item The crystal lattice $\mathbf{Cr}$ is characteristic for $\CL$,
\item Belnap's $\MM$ is characteristic for $\BM$, and
\item each of the Sugihara chains $\Su_{\Zn^*}$, $\Su_{\Zn}$, and
  $\Su_{\Zn+\Zn}$ is characteristic for $\RM$.
\end{itemize}
These lattices can be represented as algebras of subsets of relevant
model structures.  This was done for $\CL$ by two relevant model
structures, one with 45 triples of elements of $\{\mathsf{T},
\mathsf{T}^*, \mathsf{a}, \mathsf{a}^*\}$, and the other with 49, the
largest possible number of triples that can be used for this purpose
\cite[pp.\ 95--100]{MR3728341}.  Both structures produce the table on
\cite[p.\ 97]{MR3728341}.  The relevant model structure in Table
\ref{crystaltables} is
\begin{gather*}
  \<\{\ll0,\cl0,\ll1,\cl1\},\C,\conv{},\emptyset\>,\\
  \C=\{\<\a,\b,\c\>\colon\a,\b,\c\in\{\ll0,\cl0,\ll1,\cl1\},\,
  \a|'\b\supseteq\c\}.
\end{gather*} 
It has only 24 triples, and is isomorphic to the restriction of the
canonical atom structure of $\sra_{\{0,1\}}$ to the diversity atoms.
The table for $\to'$ coincides with the table on \cite[p.\
97]{MR3728341} when $\mathsf{T}=\ll1$, $\mathsf{T}^*=\cl1$,
$\mathsf{a}=\ll0$, and $\mathsf{a}^*=\cl0$.  Other numbers of triples
besides 24, 45, and 49 work, but the choice made here has the feature
that the ternary relation holds among binary relations, instead of
unstructured objects.  The ternary relation is set-theoretically
defined as ``the relativized relative product of the first two
contains the third'', and the Routley star ${}^*$ is conversion
$\conv{}$.  Similarly, $\BM$ was characterized in \cite[pp.\
100--104]{MR3728341} by a single finite relevant model structure with
13 triples of elements of $\{ \mathsf{T}, \mathsf{a}, \mathsf{a}^*
\}$, which is isomorphic to the canonical atom structure of
$\sra_{\{0\}}$.  In the notation of \SS\SS\ref{sect2}--\ref{sect3},
that structure is
\begin{gather*}
  \<\{<,>,=\},\C,\conv{},\{=\}\>\\
  \C=\{\<\a,\b,\c\>\colon\a,\b,\c\in\{<,>,=\},\,\a|\b\supseteq\c\}.
\end{gather*}
The points $\mathsf{T}$, $\mathsf{a}$, and $\mathsf{a}^*$ match up
with the binary relations $=$, $<$, and $>$ on the rationals. The
Routley-Meyer ternary relation in this case is ``the relative product
of the first two contains the third'', and the Routley star is
conversion.  The atom structure of the complete atomic proper relation
algebra $\sra_\Zn$ is
\begin{gather*}
  \<\At_\Zn,\C,\conv{},\{\Id_\Zn\}\>,\\
  \C=\{\<\a,\b,\c\>\colon\a,\b,\c\in\At_\Zn, \, \a|\b\supseteq\c\}.
\end{gather*}
This relevant model structure is characteristic for $\RM$, and its
ternary relation is the product-inclusion relation.

In the previous examples, the ternary relation of the Routley-Meyer
semantics is the product-inclusion relation, possibly relativized.
The logic $\KR$ is different. The atom structure of the canonical
extension of the free symmetric dense relation algebra on countably
many generators is a relevant model structure characteristic for
$\KR$.  The same structure can be constructed by letting the atoms be
maximal $\KR$-theories.  In both cases the atoms are not binary
relations, nor can they be represented as relations, because there are
symmetric dense relation algebras that are not representable.  For
example, there are three non-representable symmetric dense relation
algebras with four atoms, but none smaller.  The 65 symmetric relation
algebras with four atoms are numbered \alg{1}{65}--\alg{65}{65} in
\cite{MR2269199}.  The three that are non-representable and dense are
\alg{36}{65}, \alg{42}{65}, and \alg{50}{65}.  

On the other hand, many symmetric dense relation algebras \emph{are}
representable, such as Church's proper 3-atom relation algebra
$\mathfrak{Ch}$.  The atom structure of $\mathfrak{Ch}$ is a relevant
model structure verifying $\KR$ that is characteristic for a complete
decidable extension of $\KR$. Once again, its ternary relation is the
relativized product-inclusion relation and its Routley star is
conversion.  Problem \ref{prob1} asks how far this kind of analysis
can be extended. What other algebras can be represented with binary
relations?
\begin{problem}\label{prob2}
  If a relation algebra is finite, integral, possibly commutative, and
  every one of its diversity atoms is dense, transitive, and distinct
  from its converse, must that algebra be representable?
\end{problem}
\subsection*{Remarks on Problem \ref{prob2}}
There are three relation algebras with five atoms that contain the
crystal lattice.  In the numbering system of \cite{MR2269199}, they
are \alg{2}{83}, \alg{29}{83}, and \alg{43}{83} (the second,
twenty-ninth, and forty-third algebras in a list of 83 algebras in
\cite[Ch.\ 6, \SS\SS62-3]{MR2269199} whose atoms are the identity
element $\id$, plus two diversity atoms $\r$, $\s$, and their
converses $\con\r$, $\con\s$).  Algebras \alg{2}{83} and \alg{43}{83}
are commutative, but \alg{29}{83} is not commutative.  The two
commutative algebras are representable. In fact, \alg{2}{83} is
isomorphic to $\sra_{\{0,1\}}$.  What about the non-commutative
algebra \alg{29}{83}?  Is it representable? This is the smallest
particular instance of Problem \ref{prob2}.
\begin{problem}\label{prob3}
  Explore the structure of algebras in $\KRM$. Does the traditional
  axiomatic approach to $\RM$ yield a finite equational axiomatization
  for the variety generated by $\KRM$?
\end{problem}
\subsection*{Remarks on Problem \ref{prob3}}
$\sra_\I\in\KRM$ for every $\I\subseteq\Zn$.  Preliminary
investigation shows $\KRM$ has many algebras that are not linearly
ordered. What else is in $\KRM$?
\begin{problem}\label{prob4}
  Do the product-inclusion relations $\A|\B\supseteq\C$ and
  $\A|'\B\supseteq\C$ on binary relations have any bearing on the
  concepts of relevance and conditionality? Do the residuations
  $\A\to\B$ and $\A\to'\B$ have any bearing on entailment?  How do the
  De Morgan negations $\rmin\A$ and $\rmin'\A$ compare and contrast
  with the Boolean negation $\min\A$?  How are relative multiplication
  and conversion related to fusion and the Routley star?  Is $\KRM$ a
  mathematical standard example of Routley-Meyer semantics?
\end{problem}
\subsection*{Remarks on Problem \ref{prob4}}
We will comment on each of the questions in Problem \ref{prob4} in the
order they occur. We start with a ternary relation on ordered pairs.
On \cite[p.\ 599]{MR2914447}, the authors say,
\begin{quote}
  ``[I]n the semantics given by Routley and Meyer, the crucial ternary
  relation $R$ is involved in the semantics as follows: for any
  sentences $A$ and $B$ at any point $x$ in any model $M$:
  \begin{equation}\tag{R}
    \text{$x\models_M\A\to\B$ iff for all $y,z$ such that $Rxyz$, if
      $y\models_M\A$ then $z\models_M\B$.}
  \end{equation}
  \dots In order to provide a philosophically illuminating semantics
  of the relevant conditional, we need to say more about what these
  models are: what the points are, what the ternary relation $R$ is,
  and why compound sentences---in particular conditionals---are
  evaluated in the way that they are. What's more, this explication
  had better make it clear how these models relate to conditionality;
  otherwise the semantics can be fairly accused of arbitrariness, or
  of ad hocness, or of simply copying the phenomenon to be explained.
  In short, the semantics are `merely formal' and philosophically
  unilluminating---at least if we want to understand the
  \emph{meaning} of a conditional. So more is required.''
\end{quote}
Later, on \cite[p.\ 601]{MR2914447}, they suggest,
\begin{quote}
  ``What's going on (or what may be seen as such) is that our
  conditional calls for a broader perspective on our universe of
  candidate counterexamples; it calls us to recognize `pair points' in
  addition to our `old' points.''
\end{quote}
We follow their lead and use ordered pairs as pair points. Let $\R$ be
the set of triples of the form $\<\<\a,\b\>, \<\c,\a\>, \<\c,\b\>\>$,
and let $\x=\<\a,\b\>$, $\y=\<\c,\a\>$, and $\z=\<\c,\b\>$.  Then
$\R\x\y\z$ holds, so the phrase ``such that $\R\x\y\z$'' may be
deleted from \thetag{R}, and ``for all $y,z$'' can be replaced by
``for all $c$''.  We also replace $x,y,z$ with the pairs $\<\a,\b\>$,
$\<\c,\a\>$, and $\<\c,\b\>$, respectively, and, to reduce the
clutter, we drop the angle brackets, the commas, and the subscript
$M$.  The result is
\begin{align*}\tag{$\to$}
  ab&\models\A\to\B\text{ iff for all $c$, if } ca\models\A%
  \text{ then } cb\models\B.%
  \intertext{Clauses for the other connectives are}%
  \tag{$\lor$}
  ab&\models\A\lor\B\ifff ab\models\A\text{ or } ab\models\B,\\
  \tag{$\land$}
  ab&\models\A\land\B\ifff ab\models\A\text{ and } ab\models\B,\\
  \tag{$\rmin$} ab&\models\rmin\A\ifff ba\not\models\A.%
  \intertext{Commutativity, density, and transitivity, which are
    required for $\RM$, can be added in three more clauses,}%
  \tag{comm} \text{if }ab&\models\A\text{ and }bc\models\B\text{ then,
    for some $d$, }ad\models\B\text{ and }dc\models\A,\\
  \tag{dense}%
  \text{if }ab&\models\A\text{ then, for some $c$, }
  ac\models\A\text{ and } cb\models\A,\\
  \tag{trans} \text{if }ab&\models\A\text{ and } bc\models\A\text{
    then } ac\models\A.%
  \intertext{If we wish to have connectives $\neg$ for Boolean
    negation and ${}^*$ for Routley star, we add two more clauses.}
  \tag{$\neg$}
  ab&\models\neg\A\ifff ab\not\models\A,\\
  \tag{${}^*$} ab&\models\A^*\ifff ba\models\A.  \intertext{Finally,
    add a clause defining validity in $M$ (reinstated as a
    subscript),}%
  \tag{valid}%
  &\models_M\A\ifff\text{for all $a$, } aa\models_M\A.
\end{align*}
By Theorem \ref{comp}, the clauses \thetag{$\to$}, \thetag{$\lor$},
\thetag{$\land$}, \thetag{$\rmin$}, \thetag{comm}, \thetag{dense}, and
\thetag{trans} yield sound and complete semantics for $\RM$.  To get
sound and complete semantics for ``classical'' $\RM$ in the sense of
\cite{MR0363789,MR0363789a}, add connectives $\neg$ and ${}^*$ to the
language, and add clauses \thetag{$\neg$} and \thetag{${}^*$}. The
question is, what do these semantics say about conditionality, or
entailment, or Boolean and De Morgan negation? Do they say anything
about relevance?

To illustrate the connection between fusion and relative
multiplication, recall that $\circ$ is defined by $\A\circ\B =
\rmin(\A\to\rmin\B)$ \cite[p.\ 269]{MR0406756}, \cite[p.\
xxiii]{MR1223997}. When formulas are regarded as relations,
$\A\circ\B=\B|A$ by \eqref{fuse}.  Consider the triple
$\<\x,\y,\z\>=\<ab,ca,cb\>\in\R$, and note that $\{\y\}|\{\x\} =
\{ca\}|\{ab\} = \{cb\} = \{\z\}$, hence $\{\x\}\circ\{\y\}=\{\z\}$, in
conformity with the usual connection between fusion and the
Routley-Meyer ternary relation \cite[(v) p.\ 414]{MR0327515}. Of
course, the distinction between fusion $\circ$ and relative
multiplication $|$ disappears under the assumption of commutativity.

Fusion appears as an associative and sometimes commutative operation
in various algebras arising from relevance logics.  Algebraization is
mathematically illuminating, but it is open to the charge that ``\dots
algebraic characterizations \dots are merely formal, exhibiting no
connection with the intended meanings of the logical constants,''
\cite[p.\,405]{MR551278}. Algebras of subsets of relevant model
structures do interpret $\lor$ and $\land$ as union and intersection,
but the other connectives arise abstractly from the ternary relation
$R$ and the unary operation ${}^*$ according to \thetag{R} and
$x\models\rmin\A$ iff $x^*\not\models\A$. ``If the only constraint on
${}^*$ is that the resulting theory should validate the right set of
sentences, then we are indeed in the presence of merely formal model
theory,'' \cite[p.\ 410]{MR551278}, and ``\dots it is completely
obscure what meaning is given to negation in the Routley-Meyer theory
\dots,'' \cite[p.\ 408]{MR551278}.  For $\RM$, according to the
interpretation of formulas as relations, the Routley star is
conversion and the meaning of negation is $\rmin\A=\min{\conv\A}$.
Anderson and Belnap ask \cite[p.\ 345]{MR0406756}, ``How then to
interpret $\circ$?  We confess puzzlement.'' For $\RM$, the answer is
$\A\circ\B=\B|A$.  Are these answers ``merely formal''? Is the meaning
of negation ``completely obscure''? Do these answers help us
understand fusion and star?

Perhaps the semantics of $\RM$ provided by $\KRM$ is a mathematical
``standard example''.  Maybe the semantics of $\BM$ and $\CL$ provided
by $\sra_{\{0\}}$ and $\sra_{\{0,1\}}$ are also mathematical examples.
The historical difference is that the logics $\BM$ and $\CL$ were
built around $\MM$ and $\mathbf{Cr}$, unlike $\RM$, which arose
entirely through choices of axioms based on purely logical
considerations.  Nevertheless, these choices led to $\KRM$, giving
interpretations for fusion and Routley star drawn from nineteenth
century algebraic logic, rather than simply constrained so that ``the
resulting theory should validate the right set of sentences.''  Unlike
Belnap's $\MM$, the Point Algebra did not arise from relevance logic,
and would have been intensely studied even if relevance logic never
existed.  Similarly, the definition of $\KRM$ is independent of
relevance logic, in spite of having been discovered by a careful
analysis of $\RM$.

Works that may be useful for Problem \ref{prob4} include
\cite{MR0406756, MR1223997, MR2914447, BimboDunn2017, MR2459113,
  MR611476, MR675916, Burgess2005-BURNRO, MR551278, MR585425,
  MR878520, MR3297467, MR0437309, MR1786098, MR0363789, MR0363789a,
  MR693387, MR0409114, MR1352869, MR1441077, zbMATH01599882}.
\section{Concluding remarks}\label{sect12}
In the introduction \cite[\SS1]{MR2536403}, the concept of dynamic
semantics is described.
\begin{quote} 
  ``Collections of \emph{binary relations} can be viewed as a sort of
  dynamic interpretation for a logic, that is thought to describe the
  impact a sentence has on a situation via specifying a set of
  possible resulting situations. Special types of \emph{dynamic
    semantics} are those in which the binary relations constitute a
  relation algebra or a relevant family of operations.''
\end{quote}
Comparison of definitions shows that a relevant family of relations
(``operations'' was a misprint) is a direct reduct of a proper
relation algebra on a set.  The concluding remarks begin
\cite[\SS7]{MR2536403},
\begin{quote} 
  ``As the reader has surely realized by now, constructing a dynamic
  semantics with certain closure properties is not a trivial
  enterprise because of the nonrepresentability result of Lyndon
  (1950).''
\end{quote}
Of course, for non-representable relation algebras the construction of
dynamic semantics for their logics is not possible.  The first part of
this sentence was an understatement for such cases. However, no
non-representable relation algebra appears in \cite{MR2536403},
certainly not Lyndon's \cite{MR0037278}, nor is there any formula
valid in proper relation algebras that is not also a theorem of $\RR$.
For such formulas, consult \cite{MR2641636,MR2496334}.  In other
cases, such as the logic $\BM$, dynamic semantics are possible because
$\MM$ is the direct reduct of the Point Algebra, which was noted in
\cite{MR2536403}. As we have shown, there are dynamic semantics for
the logic $\CL$ and the logic of Meyer's RM84 (which has not been
axiomatized, so far as we know).

For any particular relevant model structure, it may not be readily
apparent whether it is representable as a definitional reduct of a
proper relation algebra.  When it is, its logic has a dynamic
semantics. Unless non-representability has been proved, it is not safe
to assume that dynamic semantics cannot be found; see the remarks
following Theorem \ref{main}.  Similarly, any particular formula has a
meaning if regarded as a statement about binary relations. In the
theory of relation algebras, this meaning matters.  It was the target
of the axiomatizations by McKinsey in 1940 \cite{MR2513} and Tarski in
1941 \cite{MR5280}, directed as they were at the Peirce-Schr\"oder
calculus of binary relations.

The concept of representability was present in the theory of relation
algebras from its inception, as the title of the 1948 abstract
\cite{JonssonTarski1948} makes clear, but was absent from relevance
logic until 2007. The representability of finite Sugihara matrices and
the resulting dynamic semantics for $\RM$ are only a decade old.  Even
today it is possible to construct new and unsuspected dynamic
semantics for rather old algebras and logics, as has been done here.
The 1952 J\'onsson-Tarski Representation Theorem \cite[3.10]{MR44502},
stated here as Theorem \ref{repth}, could be regarded as a successful
application of the Routley-Meyer ternary relation and the Routley
star, one that provides Routley-Meyer semantics for all relation
algebras. But the J\'onsson-Tarski Representation Theorem preceded the
introduction of Routley-Meyer semantics by two decades.  Shortly
before his death, Meyer was informed of the relational completeness
theorem for $\RM$ \cite[Theorem 6.2]{MR2641636}.  His response was an
email message that ended with ``KEEP'EM COMING''.  In this paper we
have tried to do so.
\section*{Acknowledgements}
We are very grateful to four reviewers for their careful work, their
numerous valuable suggestions, and for catching many errors in the
first version.


\begin{thebibliography}{10}

\bibitem{Allen1981}
James~F. Allen.
\newblock An interval-based representation of temporal knowledge.
\newblock In {\em Proceedings of the Seventh International Joint Conference on
  Artificial Intelligence, (IJCAI)}, pages 221--226, 1981.

\bibitem{Allen1984}
James~F. Allen.
\newblock {Towards a general theory of action and time}.
\newblock {\em Artificial Intelligence}, 23(2):123--154, July 1984.

\bibitem{Allen1983}
James~F. Allen.
\newblock Maintaining knowledge about temporal intervals.
\newblock {\em Communications of the Association for Computing Machinery},
  26(11):832--842, November 1983.

\bibitem{MR0406756}
Alan~Ross Anderson and Nuel~D. Belnap, Jr.
\newblock {\em Entailment. The logic of relevance and necessity. Vol. I}.
\newblock Princeton University Press, Princeton, N. J.-London, 1975.

\bibitem{MR1223997}
Alan~Ross Anderson, Nuel~D. Belnap, Jr., and J.~Michael Dunn.
\newblock {\em Entailment. The logic of relevance and necessity. Vol. II}.
\newblock Princeton University Press, Princeton, NJ, 1992.

\bibitem{MR1334290}
Hajnal Andr\'{e}ka and Roger~D. Maddux.
\newblock Representations for small relation algebras.
\newblock {\em Notre Dame J. Formal Logic}, 35(4):550--562, 1994.

\bibitem{MR3526497}
Arnon Avron.
\newblock R{M} and its nice properties.
\newblock In {\em J. {M}ichael {D}unn on information based logics}, volume~8 of
  {\em Outst. Contrib. Log.}, pages 15--43. Springer, [Cham], 2016.

\bibitem{MR2235660}
Silvana Badaloni and Massimiliano Giacomin.
\newblock The algebra {IA$^\text{fuz}$}: a framework for qualitative fuzzy
  temporal reasoning.
\newblock {\em Artificial Intelligence}, 170(10):872--908, 2006.

\bibitem{MR2225706}
Philippe Balbiani, Jean-Fran\c{c}ois Condotta, and G\'{e}rard Ligozat.
\newblock On the consistency problem for the {INDU} calculus.
\newblock {\em J. Appl. Log.}, 4(2):119--140, 2006.

\bibitem{MR2914447}
Jc~Beall, Ross Brady, J.~Michael Dunn, A.~P. Hazen, Edwin Mares, Robert~K.
  Meyer, Graham Priest, Greg Restall, David Ripley, John Slaney, and Richard
  Sylvan.
\newblock On the ternary relation and conditionality.
\newblock {\em J. Philos. Logic}, 41(3):595--612, 2012.

\bibitem{MR0141590}
Nuel~D. Belnap, Jr.
\newblock Entailment and relevance.
\newblock {\em J. Symbolic Logic}, 25:144--146, 1960.

\bibitem{BimboDunn2017}
K.~Bimbo and J.~M. Dunn.
\newblock The emergence of set-theoretical semantics for relevance logics
  around 1970.
\newblock {\em Proceedings of the Third Workshop, May 16--17, 2016, Edmonton,
  Canada, (IFCoLog Journal of Logics and Their Application)}, 4(3):557--589,
  2017.

\bibitem{MR2459113}
Katalin Bimb\'{o} and J.~Michael Dunn.
\newblock {\em Generalized {G}alois logics}, volume 188 of {\em CSLI Lecture
  Notes}.
\newblock CSLI Publications, Stanford, CA, 2008.

\bibitem{MR2536403}
Katalin Bimb\'{o}, J.~Michael Dunn, and Roger~D. Maddux.
\newblock Relevance logics and relation algebras.
\newblock {\em Rev. Symb. Log.}, 2(1):102--131, 2009.

\bibitem{MR3873385}
Manuel Bodirsky.
\newblock Finite relation algebras with normal representations.
\newblock In {\em Relational and algebraic methods in computer science. 17th
  international conference, RAMiCS 2018, Groningen, The Netherlands, October 29
  -- November 1, 2018. Proceedings}.

\bibitem{MR2565927}
Manuel Bodirsky and Hubie Chen.
\newblock Qualitative temporal and spatial reasoning revisited.
\newblock {\em J. Logic Comput.}, 19(6):1359--1383, 2009.

\bibitem{MR2722825}
Manuel Bodirsky and Jan K\'{a}ra.
\newblock A fast algorithm and {D}atalog inexpressibility for temporal
  reasoning.
\newblock {\em ACM Trans. Comput. Log.}, 11(3):Art. 15, 21, 2010.

\bibitem{zbMATH03928963}
Ross~T. {Brady}.
\newblock {Depth relevance of some paraconsistent logics.}
\newblock {\em {Stud. Log.}}, 43:63--73, 1984.

\bibitem{MR1909067}
Mathias Broxvall.
\newblock The point algebra for branching time revisited.
\newblock In {\em K{I} 2001: {A}dvances in artificial intelligence ({V}ienna)},
  volume 2174 of {\em Lecture Notes in Comput. Sci.}, pages 106--121. Springer,
  Berlin, 2001.

\bibitem{MR2015214}
Mathias Broxvall and Peter Jonsson.
\newblock Point algebras for temporal reasoning: algorithms and complexity.
\newblock {\em Artificial Intelligence}, 149(2):179--220, 2003.

\bibitem{MR611476}
John~P. Burgess.
\newblock Relevance: a fallacy?
\newblock {\em Notre Dame J. Formal Logic}, 22(2):97--104, 1981.

\bibitem{MR675916}
John~P. Burgess.
\newblock Common sense and ``relevance''.
\newblock {\em Notre Dame J. Formal Logic}, 24(1):41--53, 1983.

\bibitem{Burgess2005-BURNRO}
John~P. Burgess.
\newblock No requirement of relevance.
\newblock In Stewart Shapiro, editor, {\em Oxford Handbook of Philosophy of
  Mathematics and Logic}, pages 727--750. Oxford University Press, 2005.

\bibitem{MR0043763}
Louise~H. Chin and Alfred Tarski.
\newblock Distributive and modular laws in the arithmetic of relation algebras.
\newblock {\em Univ. California Publ. Math. (N.S.)}, 1:341--384, 1951.

\bibitem{MR3466223}
Jean-Fran\c{c}ois Condotta, Souhila Kaci, and Yakoub Salhi.
\newblock Optimization in temporal qualitative constraint networks.
\newblock {\em Acta Inform.}, 53(2):149--170, 2016.

\bibitem{MR551278}
B.~J. Copeland.
\newblock On when a semantics is not a semantics: some reasons for disliking
  the {R}outley-{M}eyer semantics for relevance logic.
\newblock {\em J. Philos. Logic}, 8(4):399--413, 1979.

\bibitem{MR585425}
B.~J. Copeland.
\newblock The trouble {A}nderson and {B}elnap have with relevance.
\newblock {\em Philos. Stud.}, 37(4):325--334, 1980.

\bibitem{MR878520}
B.~J. Copeland.
\newblock What is a semantics for classical negation?
\newblock {\em Mind}, 95(380):478--490, 1986.

\bibitem{zbMATH03513746}
J.~Michael {Dunn}.
\newblock {A {K}ripke-style semantics for {R}-mingle using a binary
  accessibility relation.}
\newblock {\em {Stud. Log.}}, 35:163--172, 1976.

\bibitem{MR2067967}
J.~Michael Dunn.
\newblock A representation of relation algebras using {R}outley-{M}eyer frames.
\newblock In {\em Logic, Meaning and Computation}, volume 305 of {\em Synthese
  Lib.}, pages 77--108. Kluwer Acad. Publ., Dordrecht, 2001.

\bibitem{MR3297467}
J.~Michael Dunn.
\newblock The relevance of relevance to relevance logic.
\newblock In {\em Logic and its applications}, volume 8923 of {\em Lecture
  Notes in Comput. Sci.}, pages 11--29. Springer, Heidelberg, 2015.

\bibitem{MR0437309}
Kit Fine.
\newblock Models for entailment.
\newblock {\em J. Philos. Logic}, 3:347--372, 1974.

\bibitem{MR1786098}
Dov~M. Gabbay and Heinrich Wansing, editors.
\newblock {\em What is negation?}, volume~13 of {\em Applied Logic Series}.
\newblock Kluwer Academic Publishers, Dordrecht, 1999.

\bibitem{MR1658912}
Rosella Gennari.
\newblock Temporal reasoning and constraint programming: a survey.
\newblock {\em CWI Quarterly}, 11(2-3):163--214, 1998.

\bibitem{MR2158576}
Alfonso Gerevini.
\newblock Incremental qualitative temporal reasoning: algorithms for the point
  algebra and the {ORD}-{H}orn class.
\newblock {\em Artificial Intelligence}, 166(1-2):37--80, 2005.

\bibitem{MR2790871}
Alfonso~E. Gerevini and Alessandro Saetti.
\newblock Computing the minimal relations in point-based qualitative temporal
  reasoning through metagraph closure.
\newblock {\em Artificial Intelligence}, 175(2):556--585, 2011.

\bibitem{MR3699802}
Steven Givant.
\newblock {\em Advanced topics in relation algebras---Relation Algebras. {V}ol.
  2}.
\newblock Springer, Cham, 2017.

\bibitem{MR3699801}
Steven Givant.
\newblock {\em Introduction to relation algebras---Relation Algebras. {V}ol.
  1}.
\newblock Springer, Cham, 2017.

\bibitem{MR1117874}
Steven~R. Givant.
\newblock A portrait of {A}lfred {T}arski.
\newblock {\em Math. Intelligencer}, 13(3):16--32, 1991.

\bibitem{MR1872076}
Gabrielle~Assunta Gr\"{u}n.
\newblock An efficient algorithm for the uniform maximum distance problem on a
  chain.
\newblock {\em Discrete Math. Theor. Comput. Sci.}, 4(2):323--350, 2001.

\bibitem{MR781929}
Leon Henkin, J.~Donald Monk, and Alfred Tarski.
\newblock {\em Cylindric algebras. {P}art {I}}, volume~64 of {\em Studies in
  Logic and the Foundations of Mathematics}.
\newblock North-Holland Publishing Co., Amsterdam, 1985.
\newblock With an introductory chapter: General theory of algebras, Reprint of
  the 1971 original.

\bibitem{MR0124250}
Leon Henkin and Alfred Tarski.
\newblock Cylindric algebras.
\newblock In {\em Proc. {S}ympos. {P}ure {M}ath., {V}ol. {II}}, pages 83--113.
  American Mathematical Society, Providence, R.I., 1961.

\bibitem{MR1935083}
Robin Hirsch and Ian Hodkinson.
\newblock {\em Relation algebras by games}, volume 147 of {\em Studies in Logic
  and the Foundations of Mathematics}.
\newblock North-Holland Publishing Co., Amsterdam, 2002.
\newblock With a foreword by Wilfrid Hodges.

\bibitem{MR3944722}
Robin Hirsch, Marcel Jackson, and Tomasz Kowalski.
\newblock Algebraic foundations for qualitative calculi and networks.
\newblock {\em Theoret. Comput. Sci.}, 768:99--116, 2019.

\bibitem{MR3089976}
Yoichi Iwata and Yuichi Yoshida.
\newblock Exact and approximation algorithms for the maximum constraint
  satisfaction problem over the point algebra.
\newblock In {\em 30th {I}nternational {S}ymposium on {T}heoretical {A}spects
  of {C}omputer {S}cience}, volume~20 of {\em LIPIcs. Leibniz Int. Proc.
  Inform.}, pages 127--138. Schloss Dagstuhl. Leibniz-Zent. Inform., Wadern,
  2013.

\bibitem{JonssonTarski1948}
Bjarni J{\'o}nsson and Alfred Tarski.
\newblock Representation problems for relation algebras.
\newblock {\em Bulletin of the American Mathematical Society}, 54:80 and 1192,
  1948.
\newblock Abstract 89.

\bibitem{MR44502}
Bjarni J\'{o}nsson and Alfred Tarski.
\newblock Boolean algebras with operators. {I}.
\newblock {\em Amer. J. Math.}, 73:891--939, 1951.

\bibitem{MR45086}
Bjarni J\'{o}nsson and Alfred Tarski.
\newblock Boolean algebras with operators. {II}.
\newblock {\em Amer. J. Math.}, 74:127--162, 1952.

\bibitem{MR1369206}
Peter~B. Ladkin and Roger~D. Maddux.
\newblock On binary constraint problems.
\newblock {\em J. Assoc. Comput. Mach.}, 41(3):435--469, 1994.

\bibitem{MR0037278}
Roger~C. Lyndon.
\newblock The representation of relational algebras.
\newblock {\em Ann. of Math. (2)}, 51:707--729, 1950.

\bibitem{MR79570}
Roger~C. Lyndon.
\newblock The representation of relation algebras. {II}.
\newblock {\em Ann. of Math. (2)}, 63:294--307, 1956.

\bibitem{MR662049}
Roger~D. Maddux.
\newblock Some varieties containing relation algebras.
\newblock {\em Trans. Amer. Math. Soc.}, 272(2):501--526, 1982.

\bibitem{MR2269199}
Roger~D. Maddux.
\newblock {\em Relation Algebras}, volume 150 of {\em Studies in Logic and the
  Foundations of Mathematics}.
\newblock Elsevier B. V., Amsterdam, 2006.

\bibitem{Maddux2007}
Roger~D. Maddux.
\newblock Relevance logic and the calculus of relations (abstract).
\newblock In {\em International Conference on Order, Algebra, and Logics,
  Vanderbilt University, June 13, 2007}, pages 1--3, 2007.

\bibitem{MR2641636}
Roger~D. Maddux.
\newblock Relevance logic and the calculus of relations.
\newblock {\em Rev. Symb. Log.}, 3(1):41--70, 2010.

\bibitem{MR2616328}
Ralph Nelson~Whitfield McKenzie.
\newblock {\em T{HE} {REPRESENTATION} {OF} {RELATION} {ALGEBRAS}}.
\newblock ProQuest LLC, Ann Arbor, MI, 1966.
\newblock Thesis (Ph.D.)--University of Colorado at Boulder.

\bibitem{MR2513}
J.~C.~C. McKinsey.
\newblock Postulates for the calculus of binary relations.
\newblock {\em J. Symbolic Logic}, 5:85--97, 1940.

\bibitem{MR6334}
J.~C.~C. McKinsey.
\newblock A solution of the decision problem for the {L}ewis systems {S}2 and
  {S}4, with an application to topology.
\newblock {\em J. Symbolic Logic}, 6:117--134, 1941.

\bibitem{MR0327515}
Robert~K. Meyer and Richard Routley.
\newblock Algebraic analysis of entailment. {I}.
\newblock {\em Logique et Analyse (N.S.)}, 15:407--428, 1972.

\bibitem{MR0363789}
Robert~K. Meyer and Richard Routley.
\newblock Classical relevant logics. {I}.
\newblock {\em Studia Logica}, 32:51--68, 1973.

\bibitem{MR0363789a}
Robert~K. Meyer and Richard Routley.
\newblock Classical relevant logics. {II}.
\newblock {\em Studia Logica}, 33:183--194, 1974.

\bibitem{MR2496334}
Szabolcs Mikul\'{a}s.
\newblock Algebras of relations and relevance logic.
\newblock {\em J. Logic Comput.}, 19(2):305--321, 2009.

\bibitem{MR1927631}
I.~Navarrete, A.~Sattar, R.~Wetprasit, and R.~Marin.
\newblock On point-duration networks for temporal reasoning.
\newblock {\em Artificial Intelligence}, 140(1-2):39--70, 2002.

\bibitem{MR693387}
R.~Routley, V.~Routley, R.~K. Meyer, and E.~P. Martin.
\newblock On the philosophical bases of relevant logic semantics.
\newblock {\em J. Non-Classical Logic}, 1(1):71--105, 1982.

\bibitem{MR544617}
Richard Routley.
\newblock Alternative semantics for quantified first degree relevant logic.
\newblock {\em Studia Logica}, 38(2):211--231, 1979.

\bibitem{MR0409114}
Richard Routley and Robert~K. Meyer.
\newblock The semantics of entailment. {I}.
\newblock In {\em Truth, syntax and modality (Proc. Conf. Alternative
  Semantics, Temple Univ., Philadelphia, Pa., 1970)}, volume~68 of {\em Studies
  in Logic and the Foundations of Mathematics}, pages 199--243. North-Holland,
  Amsterdam, 1973.

\bibitem{MR728950}
Richard Routley, Val Plumwood, Robert~K. Meyer, and Ross~T. Brady.
\newblock {\em Relevant logics and their rivals. {P}art {I}}.
\newblock Ridgeview Publishing Co., Atascadero, CA, 1982.
\newblock The basic philosophical and semantical theory.

\bibitem{Schechter2005}
Eric Schechter.
\newblock {\em Classical and nonclassical logics}.
\newblock Princeton University Press, Princeton, NJ, 2005.

\bibitem{Sugihara1955}
Takeo Sugihara.
\newblock Strict implication free from implicational paradoxes.
\newblock {\em Memoirs of the Faculty of Liberal Arts, Fukei University},
  Series I:55--59, 1955.

\bibitem{MR3728341}
Richard Sylvan, Robert Meyer, Val Plumwood, and Ross Brady.
\newblock {\em Relevant logics and their rivals. {V}ol. {II}}, volume~59 of
  {\em Western Philosophy Series}.
\newblock Ashgate Publishing Limited, Aldershot, 2003.
\newblock A continuation of the work of Richard Sylvan, Robert Meyer, Val
  Plumwood and Ross Brady, Edited by Brady.

\bibitem{zbMATH02525153}
Alfred {Tarski}.
\newblock {Grundz\"uge des Systemenkalk\"uls. I.}
\newblock {\em {Fundam. Math.}}, 25:503--526, 1935.

\bibitem{MR5280}
Alfred Tarski.
\newblock On the calculus of relations.
\newblock {\em J. Symbolic Logic}, 6:73--89, 1941.

\bibitem{MR736686}
Alfred Tarski.
\newblock {\em Logic, semantics, metamathematics}.
\newblock Hackett Publishing Co., Indianapolis, IN, second edition, 1983.
\newblock Papers from 1923 to 1938, Translated by J. H. Woodger, Edited and
  with an introduction by John Corcoran.

\bibitem{vanBenthem1984}
Johan van Benthem.
\newblock Review: On when a semantics is not a semantics: Some reasons for
  disliking the {R}outley-{M}eyer semantics for relevance logic, by {B}. {J}.
  {C}opeland.
\newblock {\em The Journal of Symbolic Logic}, 49(3 Sep., 1984):994--995.

\bibitem{MR1352869}
Johan van Benthem.
\newblock {\em Language in action}.
\newblock MIT Press, Cambridge, MA; copublished with North-Holland, Amsterdam,
  1995.
\newblock Categories, lambdas and dynamic logic, Revised reprint of the 1991
  original.

\bibitem{MR1441077}
Johan van Benthem.
\newblock {\em Exploring logical dynamics}.
\newblock Studies in Logic, Language and Information. CSLI Publications,
  Stanford, CA; FoLLI: European Association for Logic, Language and
  Information, Amsterdam, 1996.

\bibitem{zbMATH01599882}
Heinrich {Wansing}, editor.
\newblock {\em {Negation. A notion in focus. Proceedings volume of an
  interdisciplinary workshop held during the conference Analyomen 2 of the
  Gesellschaft f\"ur Analytische Philosophie (GAP), Leipzig, Germany, September
  7--10, 1994.}}
\newblock Berlin: De Gruyter, 1996.

\end{thebibliography}
\end{document}